\documentclass[10pt,a4paper]{article}

\usepackage{amssymb,amsmath,amsthm,amscd,enumerate}
\usepackage[all]{xy}
\usepackage[french]{babel}

\advance\oddsidemargin -1,5cm
\advance\evensidemargin -2cm
\newcommand{\n}{\noindent}
\newcommand{\OO}{\mathcal{O}_{X}}
\newcommand{\sy}{{\mathcal Sym}}
\newcommand{\W}{\mathsf{\Lambda}}
\newcommand{\A}{{\mathcal A}}
\newcommand{\E}{{\mathcal E}}
\newcommand{\I}{{\mathcal I}}
\newcommand{\LL}{{\mathcal L}}
\newcommand{\oo}{{\mathcal O}}
\newcommand{\rr}{\longrightarrow}
\newcommand{\pp}{{\bf P}}

\begin{document}
\begin{center}
{\Large {\bf Note on double coverings and binary quadratic forms }}
\vspace{0,5cm}

Daniel Ferrand
\end{center}
\vspace{5mm} 

\begin{quotation}
{\sc Abstract}\; Let $E$ be a rank two vector bundle on a scheme $X$. The following three structures are shown to be equivalent :

\n a) A quadratic map $q : E \rr L$, with values in an invertible $\OO$-module $L$ ($q$ is assumed to be everywhere a non zero map).

\n b) A double covering $f: Y \rightarrow X$ endowed with an invertible $\mathcal{O}_{Y}$-module $\E$, plus an isomorphism $f_{\star}\E \simeq E$.

\n c) An effective Cartier divisor on     ${\bf P}_{X}(E)$, of degree two over $X$.

The passages from one of these points of view to another, although likely to specialists, deserved to be  carefully settled in their greatest generality: we only need that 2 is invertible on $X$. The passage from b) to a) puts the "norm form" in front. In the last two paragraphs the base scheme $X$ is the projective space ${\bf P}_{n}$ ; we prove that for any double covering $Y \rightarrow {\bf P}_{n}$, the homomorphism on Picard groups it induces is an isomorphism if $n \geq 3$; we finally apply this result to quadratic forms on rank two vector bundles on ${\bf P}_{n}$.

\bigskip

\n {\sc MSC 2010} : 11E16; 14J60; 32L10.
\end{quotation}
\bigskip

\noindent {\bf {\large Introduction}}
\medskip

This note deals with rank two vector bundles on a scheme $X$. In a first part, we characterize those bundles which are the direct image $f_{\star}(\E)$ of an invertible  module $\E$ defined on some double covering $f : Y \rr X$ (All the coverings we consider are finite and flat, and usually of rank 2 ; no assumptions are made concerning their ramification, nor the possible singularities of $X$ or $Y$ ; but we always suppose that 2 is everywhere invertible). The point is that a direct image of this kind comes equipped with the extra structure given by the "norm" ({\it i.e.} the wedge 2): namely one has the map
$$
\nu_{\E}: f_{\star}\E \; \rr \; {\sf N}_{Y/X}(\E) = {\mathcal Hom}(\W^2f_{\star}\oo_{Y}, \W^2f_{\star}\E),\leqno{(\star)}
$$
defined by
$$
x \; \longmapsto \; \wedge^2(\alpha \mapsto \alpha x)
$$
It is a quadratic map with values in the invertible $\OO$-module ${\sf N}_{Y/X}(\E)$. Conversely, given a locally free $\OO$-module $E$ of rank 2, and  a \emph{surjective} linear map
$$
\varphi : \sy^2(E) \rr L,
$$
where $L$ is invertible, we construct a double covering $f :Y \rightarrow X$ together with an invertible $\oo_{Y}$-module $\E$ such that the quadratic map $q : E \rightarrow L$ associated with $\varphi$, is isomorphic to the "norm" written above; in particular, $f_{\star}\E$ is isomorphic to $E$. This construction goes as follows: letting 
$N = {\mathcal Hom}(L, \W^2E)$, we extract from $\varphi$, in the usual way, a linear map \, $u: N\otimes E \rightarrow E$; its wedge 2 gives a linear map $\mu: N^{\otimes 2} \rightarrow \OO$, and thus an $\OO$-algebra structure on $\A = \OO \oplus N$; moreover, the map $u$ extends to an $\A$-module structure on $E$; then we take for $f$ the double covering ${\mathcal Spec}(\A) \rr X$. This covering is \'etale if and only if $q$ is non degenerate.
\medskip

Establishing such a correspondence  required a precise description/construction of double coverings; so the first paragraph gathers results on them; the only - but essential - assumption we need is that 2 is invertible. In the context of commutative rings we can skim over this description as follows: a "quadratic algebra" $R \rightarrow S$ ({\it i.e.} as an $R$-module, $S$ is  projective of rank 2) has a direct sum decomposition $S = R \oplus N$, where $N$ is  the kernel of the trace map $ N = {\rm Ker}({\rm Tr}_{S/R})$; the multiplication is given by the map $\mu: N^{\otimes 2} \rightarrow R$ defined by $\mu(\alpha^{\otimes 2}) = - {\rm norm}_{S/R}(\alpha)$. In some sense, $S$ thus  "represents" the square root of $\mu$, finding the usual intuition back. But, expliciting this map $\mu$, may sometimes be difficult.
\n By way of example, the \S 3 gives a comprehensive description of the invertible module $N$ and of the mutiplication $\mu$ for the double covering 
$$
({\bf P}_{1})^n / \mathfrak{A}_{n}\; \longrightarrow \; ({\bf P}_{1})^n / \mathfrak{S}_{n} = {\bf P}_{n}
$$
We find that $N \simeq \oo_{{\bf P}_{n}}(1-n)$; roughly speaking, this covering represents (a globalization of) extracting square root of the discriminant of the generic polynomial of degree $n$.
\medskip

\n The correspondence alluded to above is detailed in the \S 4 and 5, and it is expanded as an equivalence between suitable categories. The main point  is the following (it is expressed here for rings) : let $S = R \oplus N$ be a quadratic algebra, and let $E$ be an invertible $S$-module. Denote by $\nu_{E} : E \rr L = {\rm Hom}(N, \W^2E)$ the "norm" as above, where we use the isomorphism $N \simeq \W^2S, \quad \alpha \mapsto 1\wedge \alpha$. Then, for $\alpha \in N$ and $x \in E$, one has
$$
x \wedge \alpha x \; = \; \nu_{E}(x)(\alpha).
$$

This formula is the key for the correspondence: since the rank of the $R$-module $E$ is  2,  once you know $x \wedge \alpha x$ for all $x$ then you also know $\alpha x$, and then the $R$-linear endomorphism $\alpha_{E} : x \mapsto \alpha x$; thus, the $S$-module structure on $E$ - that is the $\alpha_{E}$'s -  and the "norm" $\nu_{E}$ determine each other; more, the invertible module $N$ itself is given by the quadratic map $E \rightarrow L$, since $N \simeq {\rm Hom}(L, \W^2E)$, and finally, the multiplication $\mu : N^{\otimes 2} \rightarrow R$ is also computable from  $\nu_{E}$.
\medskip

\n In the \S 6 and 7, we adopt a more geometrical point of view. We consider a rank 2 vector bundle $E$ on $X$ and the projection $p : P = {\bf P}_{X}(E) \rr X$. An effective Cartier divisor $D \subset P$ such that $D \rightarrow X$ is a double covering is given by a section $\gamma :\oo_{P}\rr \oo_{P}(2)$. We determine the $\OO$-algebra structure of $p_{\star}(\oo_{D})$ from $\gamma$, via the duality between the sections of $p_{\star}(\oo_{P}(2)) = \sy^2(E)$ and the linear maps $\sy^2(E) \rightarrow \OO $.
\medskip

That is obviously reminiscent of the classical result of Schwarzenberger which states that \emph{any} rank 2 vector bundle on a projective \emph{surface} is indeed the direct image of an invertible sheaf on a double covering. In view of the above correspondence, this result is "explained" by the fact that, in dimension 2, one can always find a surjective map $\sy^2(E) \rightarrow L$ for a suitable invertible module $L$. One of the motivations in writing this note was to investigate what remains true on a scheme of greater dimension. Unfortunately, already for ${\bf P}_{n}$, with $n \geq 3$, all these direct image bundles are decomposable. In fact, we prove in the \S 8 that for \emph{any}  double covering $f : Y \rr {\bf P}_{n}$, with $n \geq 3$, the map
$$
f^{\star} : {\rm Pic}({\bf P}_{n}) \rr {\rm Pic}(Y)
$$
is an isomorphism (\emph{any} means that the scheme $Y$ may be very singular).
\medskip

The \S 9 contains an application to rank 2 bundles $E$ on ${\bf P}_{n}$. If $E$ is indecomposable and is equipped with a linear map $ \sy^2(E) \rr \oo(r)$, which is surjective in codimension $\leq 2$, then one has $r > c_{1}(E)$.
\bigskip

After the reading of a first draft of this note, Manuel {\sc Ojanguren} drew my attention to the paper [Knes], and he urged me to follow the idea of {\sc Kneser} in using  the Clifford algebras instead of assuming that 2 is invertible. In fact, given a quadratic map $q: E \rightarrow \OO$ (with $E$ of rank 2), the Clifford algebra $C$ of $q$ breaks down in the direct sum $C = C^+ \oplus C^-$, where $C^+$ is a quadratic $\OO$-algebra, and $C^-$, which is isomorphic to $E$, is a $C^+$ module. Thus, in the case where $L = \OO$, our \S 4 is contained in the Clifford algebras theory. But we definitely need to consider forms with values in an invertible module $L$ different from $\OO$; thanks to {\sc M.-A. Knus} a Clifford theory does exist also in that case,  but it seems not to be as popular as the usual one, mainly because we don't then get \emph{algebras} but a slightly more complicated structure. So, after some attempts, I gave up rewriting this text and I maintain the assumption that 2 is invertible. I wanted to keep this note as elementary as possible (except, perhaps, in the last two paragraphs), because it is intended to be nothing but preliminaries to another works on vector bundles.
\vspace{1cm}

\n {\bf 0.}\; Conventions et rappels

\n {\bf 1.}\; Rev\^etements de rang deux

\n {\bf 2.}\; Exemples

\n {\bf 3.}\; Le rev\^etement double $({\bf P}_{1})^n / \mathfrak{A}_{n}\; \longrightarrow \; {\bf P}_{n}$

\n {\bf 4.}\; Rev\^etement  double attach\'e \`a une forme quadratique

\n {\bf 5.}\; Forme quadratique sur l'image directe d'un inversible

\n {\bf 6.}\;  Formes quadratiques et polyn\^omes homog\`enes de degr\'e deux

\n {\bf 7.}\;  Diviseurs de degr\'e deux sur les fibr\'es en droites

\n {\bf 8.}\;  Groupes de Picard d'un rev\^etement de l'espace projectif

\n {\bf 9.}\;  Application aux fibr\'es de rang deux sur les espaces projectifs

\newpage
%%%%%%%%%%%%%%%%%%%%%%%%%%

\n {\bf 0. \; Conventions et rappels}
\medskip

\begin{tabular}{|c|}
\hline 
On  suppose que 2 est inversible dans tous les anneaux qui interviennent dans ce texte.\\
\hline
\end{tabular}
\bigskip

 On nommera \emph{rev\^etement} tout morphisme fini localement libre de sch\'emas, \'eventuellement ramifi\'e, voire radiciel.
\vspace{1cm}

\n {\bf 0.1.}\; Formes quadratiques

\medskip

Soient $E$ et $L$ deux modules sur un anneau commutatif $R$. Rappelons qu'on nomme  \emph{application quadratique}  une application 
$$
q : E \rr L
$$
ayant les deux propri\'et\'es suivantes :

\n {\it i}) Pour $a \in R$ et $x \in E$, on a $q(ax) = a^2q(x)$ ;

\n {\it ii}) L'application $E \times E \rr L, \quad (x, y) \longmapsto q(x+y) - q(x) - q(y)$\; est bilin\'eaire.
\medskip

\n Une telle application est la restriction \`a $E$ d'une \emph{loi polyn\^ome homog\`ene de degr\'e 2} ([Roby] p.236) ; ainsi, il y a une bijection entre l'ensemble des  applications  quadratiques $q : E \rr L$ et l'ensemble des applications lin\'eaires   $\varphi : \Gamma^2(E) \rr L$, o\`u  $\Gamma^2(E)$ d\'esigne  le module des \emph{carr\'es divis\'es}([Roby] p. 266) ; cette correspondance est d\'efinie par
$q(x) = \varphi(\gamma^{2}(x))$.
\medskip

\n Par ailleurs, 2 \'etant ici suppos\'e inversible, l'application canonique ${\rm Sym}^2(E) \rr \Gamma^2(E)$ est bijective, d'inverse d\'efini par $\gamma^2(x) \mapsto \frac{1}{2}x^2$ ($x^2$ d\'esigne le carr\'e dans   ${\rm Sym}^2$). Finalement, se donner une application quadratique $q : E \rr L$ revient  donc \`a se donner une application lin\'eaire $\varphi : {\rm Sym}^2(E) \rr L$ ; le passage de l'une \`a l'autre se voit sur les relations
$$
\varphi(xy) = q(x+y) - q(x) - q(y), \qquad q(x) = \frac{1}{2}\varphi(x^2)
$$
\bigskip

\n {\bf 0.2.}\; Pour tout $R$-module $E$, le module ${\rm Sym}^2(E)$ est engendr\'e par les carr\'es puisque, pour $x, y \in E$, on a $ xy = \frac{1}{2}((x+y)^2 - x^2 -y^2)$
\bigskip

\n {\bf 0.3.}\; Soit $L$ un $R$-module inversible, de sorte que ${\rm Sym}^2(L) = L^{\otimes 2}$. Soit $\mu : L^{\otimes 2} \rr R$ une application lin\'eaire. Alors, pour $\alpha, \beta \in L$, on a
$$
\mu(\alpha \otimes \beta) = \frac{1}{2}(\mu((\alpha+\beta)^{\otimes 2}) - \mu(\alpha^{\otimes 2}) - \mu(\beta^{\otimes 2})).
$$
\medskip

\n {\bf 0.4.}\; Partout non nulle
\medskip

Si $E$ et $F$ sont deux modules sur un anneau $R$, on dit qu'une application lin\'eaire $u : E \rr F$ est \emph{partout non nulle} si pour tout id\'eal premier $\mathfrak{p}$ de $R$, l'application $\kappa(\mathfrak{p})$-lin\'eaire $\kappa(\mathfrak{p})\otimes_{R} E \rr \kappa(\mathfrak{p})\otimes_{R} F$ est non nulle.  Si $F$ est un $R$-module inversible, il revient au m\^eme de dire que $u$ est partout non nulle, ou qu'elle est surjective.

\n Cette d\'efinition, qui concerne a priori des applications lin\'eaires, s'\'etend sans changement  aux \emph{lois polyn\^omes}, et en particulier aux applications quadratiques $q : E \rr L$ ({\bf 0.1}) ; si $L$ est un module inversible, $q$ est partout non nulle si et seulement si l'application lin\'eaire associ\'ee $\varphi : {\rm Sym}^2(E) \rr L$ est surjective. On dit alors parfois que $q$ est \emph{primitive}.

\n Cette d\'efinition s'\'etend aussi aux applications llin\'eaires entre modules quasi-coh\'erents sur un sch\'ema.

\vspace{1cm}

%%%%%%%%%%%%%%%%%%%%%%%%%

\noindent {\bf 1.\; Rev\^etements de rang deux}
\bigskip

On se propose de d\'ecrire ici les rev\^etements doubles, c'est-\`a-dire les morphismes de sch\'emas
$$
f\, :\, Y \longrightarrow \; X
$$
qui sont finis, localement libres de rang deux.

\noindent On ne fait aucune hypoth\`ese de r\'egularit\'e sur les sch\'emas $X$ et $Y$, ni sur le morphisme $f$, mais on suppose toujours que 2 est partout inversible.
\bigskip

\noindent {\bf 1.1.}\; Posons 
$$
\A \; = \; f_{\star}({\mathcal O}_{Y})
$$
C'est une $\OO$-alg\`ebre localement libre de rang 2. \`A ce titre, elle poss\`ede des applications trace et norme ; la trace
$${\rm Tr} : \A\; \rightarrow \; \OO,
$$ 
est surjective puisque ${\rm Tr}(1) = 2$ et que  2 est suppos\'e inversible. Introduisons le $\OO$-module inversible
$$
N \; = \; {\rm Ker}({\rm Tr}) \; \subset \A
$$
On a donc une d\'ecomposition en somme directe de $\OO$-modules
$$
\A \; = \; \OO \oplus N. \leqno{(1.1.1)}
$$
Utilisant cette d\'ecomposition, la multiplication dans $\A$ s'exprime simplement par une application lin\'eaire $\mu : N^{\otimes 2} \rr \OO$, de la fa\c{c}on suivante (au dessus d'un ouvert affine) : pour $a, b$ des sections  de $\OO$, et $x, y$ des sections  de $N$, on a
$$
(a+x)(b+y) \; =\; ab + \mu(x\otimes y) \;  + \;  ay+bx.
$$
 En effet, pour une section  $x \in N$, le th\'eor\`eme de Hamilton-Cayley et la relation ${\rm Tr}(x) = 0$ donnent
$$
x^2 \; = \; - {\rm norm}(x).
$$
Par suite, le produit dans $\A$ de deux \'el\'ements $x, y \in N$, est \'egal \`a
$$
xy = \frac{1}{2}((x+y)^2 - x^2 - y^2) = - \frac{1}{2}({\rm norm}(x+y) -{\rm norm}(x) - {\rm norm}(y))
$$
Le membre de droite d\'efinit  bien une forme bilin\'eaire sym\'etrique sur $N$, d'o\`u l'application
$$
\mu : N^{\otimes 2} \rr \OO. \leqno{(1.1.2)}
$$
Elle est reli\'ee \`a la norme par la formule 
$$
\mu(x\otimes y) \; = \;  - \frac{1}{2}({\rm norm}(x+y) -{\rm norm}(x) - {\rm norm}(y))
$$
\bigskip

\n {\bf 1.2.  D\'efinition }\; {\it Pour un rev\^etement double $f : Y \rr X$, le $\OO$-module inversible $$ N = {\rm Ker}({\rm Tr} : f_{\star}({\mathcal O}_{Y}) \rr \OO)$$  sera dit \emph{associ\'e \`a}  $f$, et l'application $ \mu : N^{\otimes 2} \rightarrow \OO$ sera appel\'ee \emph{ la multiplication de $Y$}.}
\bigskip

\noindent Tout couple $(N, \mu)$ form\'e d'un $\OO$-module inversible et d'une application lin\'eaire $ \mu : N^{\otimes 2} \rightarrow \OO$ d\'etermine une structure de $\OO$-alg\`ebre sur  $\OO \oplus N$, et, par suite, un morphisme ${\mathcal Spec}(\OO \oplus N) \rightarrow X$ qui est un rev\^etement de rang deux.. 
\medskip

\noindent  L'automorphisme canonique  $\sigma$ de $\A$ ([Bour] A III.13 Prop. 2) est donn\'e par $\sigma(a+x) = a-x.$ On a $\A^{\sigma} = \OO$.
\bigskip

\noindent {\bf 1.3.}\; Diramation
\medskip

\n D\'esignons par ${\mathcal I} \subset \OO$ l'image de l'application $\mu$ ; c'est un id\'eal, \'eventuellement nul, de $\OO$. Le sous-sch\'ema ferm\'e  $\Delta \subset X$ d\'efini par ${\mathcal I}$ s'appelle selon l'\'epoque, le pays ou l'auteur, le lieu de \emph{diramation}, ou de \emph{ramification}, ou de \emph{branchement} de $f$. J'utiliserai le premier terme.
\medskip

\n On a donc une suite exacte
$$
N^{\otimes 2} \; \stackrel{\mu}{\longrightarrow} \; \OO \; \longrightarrow \; {\mathcal O}_{\Delta}
\; \longrightarrow 0
$$

 L'id\'eal de $\A$ engendr\'e par $N$ est \'egal \`a ${\mathcal I} \oplus N$.  L'homomorphisme $\OO / \I  \longrightarrow \A/ \I +N$ est un isomorphisme. D'ailleurs, l'image r\'eciproque $f^{-1}(\Delta) \subset Y$ du lieu de diramation  a pour alg\`ebre $\A / \I\A = {\mathcal O}_{\Delta} \oplus N/ \I N$, le second facteur \'etant id\'eal de carr\'e nul ; en particulier, le rev\^etement  $f^{-1}(\Delta) \rightarrow \Delta$ poss\`ede une section canonique qui permet de voir $\Delta$ aussi comme un ferm\'e de $Y$.
\medskip
 
 \n La suite exacte ci-dessus montre que $\Delta$ est un diviseur sur $X$ si et seulement si $\mu$ est injective.
 Si $\mu$ est injective, alors $\I \oplus N$ est un id\'eal inversible de $\A$.
\bigskip

\n {\bf 1.4.}\; Ramification
\medskip

\n On vient de voir que $N/\I N$ peut \^etre vu comme un $\A$-module, annul\'e par $\I\A = \I \oplus \I N$.
\medskip

\n {\bf 1.4.1. Lemme }\; {\it La diff\'erentielle universelle $d_{Y/X} : \oo_{Y} \rr \Omega^1_{Y/X}$ est l'application $\, \OO \oplus N \rr N/\I N$, donn\'ee par  $ a+x \longmapsto {\rm cl}(x)$. En particulier, le lieu de diramation de $f$ est le support (sch\'ematique) de $f_{\star}(\Omega^1_{Y/X})$ ; en d'autres termes, on a $\I = {\rm Ann}_{\OO}(f_{\star}(\Omega^1_{Y/X}))$.}
{\it Le morphisme $f : Y \rr X$  est \'etale si et seulement si la multiplication $\mu : N^{\otimes 2} \longrightarrow \ \OO$ est surjective.}
\medskip

\n
Une $\OO$-d\'erivation  $D : \A= \OO \oplus N \rr E$, \`a valeurs dans un $\A$-module $E$, est en particulier  une application $\OO$-lin\'eaire $N \rightarrow E$, et il faut voir qu'elle est nulle sur $\I N$ ;  consid\'erons le produit, dans $\A$ de trois \'el\'ements $x, y, z \in N$. Comme les produits deux \`a deux de ces \'el\'ements sont dans $\OO$, on a 
$$
2 D(xyz) = D(xyz) + D(xzy) = xyD(z) + xzD(y) = x(yD(z)+zD(y)) = xD(yz) = 0.
$$
Comme $2$ est inversible, $D$ est nulle sur les produits de trois \'el\'ements de $N$, donc sur $\I N$.

\n Il reste \`a v\'erifier que l'application indiqu\'ee  $a+x \longmapsto {\rm cl}(x)$  est une d\'erivation, ce qui est \'evident compte tenu de la structure de $\A$-module sur $N/\I N$ qui  est donn\'ee par
$$
(a+x)\, , \, {\rm cl}(y) \longmapsto {\rm cl}(ay).
$$

\bigskip

\noindent {\bf 1.5.}\; Propri\'et\'e universelle.
\bigskip

\noindent Soit $f : Y \rightarrow X$ un rev\^etement double, $N$ son module inversible associ\'e et $\mu : N^{\otimes 2} \rightarrow \OO$ sa multiplication. Alors $Y$ repr\'esente le foncteur ${\sf F}$, en les sch\'emas $g : Z \rightarrow X$, d\'efini par
$$
{\sf F}(Z) =  \{ \omega : g^{\star}(N) \rightarrow {\mathcal O}_{Z}, \; \omega^{\otimes 2} = g^{\star}(\mu) \}
$$
En effet un \'el\'ement $\omega \in {\sf F}(Z)$ donne, par adjonction, une application $\OO$-lin\'eaire $\omega' : N \rightarrow g_{\star}(\oo_{Z})$ ; la condition  $\omega^{\otimes 2} = g^{\star}(\mu) $ implique que $\omega'^{\otimes 2}$ se factorise en 
$$
N^{\otimes 2} \stackrel{\mu}{\longrightarrow} \OO \rr g_{\star}(\oo_{Z})
$$
On a donc un morphisme de $\OO$-alg\`ebres $\OO \oplus N \rr g_{\star}(\oo_{Z})$, d'o\`u, finalement, un morphisme de sch\'emas $Z \rightarrow {\mathcal Spec}(\OO \oplus N) =Y$. On v\'erifie imm\'ediatement que l'application ainsi construite ${\sf F}(Z) \rr {\rm Hom}_{X}(Z, Y)$  est bijective.
\bigskip

On peut  adopter un point de vue plus syst\'ematique, en introduisant le fibr\'e vectoriel ${\bf V}_{X}(N)$ qui repr\'esente les formes lin\'eaires $N \rightarrow \oo$ ([EGA I], (9.4.8) ). L'application $\omega : f^{\star}(N) \rightarrow \mathcal{O}_{Y}$ conduit \`a un morphisme
$$
Y \; \rr \; {\bf V}_{X}(N)
$$
qui est une immersion ferm\'ee puisque l'application ${\mathcal Sym}_{\OO}(N) \rr \A = \OO \oplus N$ est surjective ; par ailleurs, l'application \og \'el\'evation au carr\'e\fg\; se traduit par un morphisme 
$$
{\bf V}_{X}(N)\; \rr \; {\bf V}_{X}(N^{\otimes 2})
$$
En termes de faisceaux d'alg\`ebres, ce morphisme est associ\'e \`a l'inclusion
$
\sy(N^{\otimes 2}) \; \subset \; \sy(N)
$ ; on a \'evidemment  $\sy(N) = \sy(N^{\otimes 2}) \oplus (N\otimes \sy(N^{\otimes 2}))$ ; ce morphisme est donc le rev\^etement double \emph{universel} associ\'e au $\sy(N^{\otimes 2})$-module inversible $N\otimes \sy(N^{\otimes 2})$.

\n La multiplication  $\mu : N^{\otimes 2} \rightarrow \OO$ conduit \`a  un morphisme de $X$-sch\'emas  $ \tilde{\mu} : X \rr {\bf V}_{X}(N^{\otimes 2})$, et le carr\'e suivant est cart\'esien:
$$
\begin{CD}
Y @>>> {\bf V}_{X}(N)\\
@VfVV  @VVV\\
X @>>\tilde{\mu} > {\bf V}_{X}(N^{\otimes 2})
\end{CD}
$$

\bigskip

\n {\bf 1.6.}\; Morphismes
\bigskip

\n {\bf Proposition}\;{\it Soient $f : Y \rightarrow X$ et $f' : Y' \rightarrow X$ deux rev\^etements doubles de $X$ ; notons $N$ et $N'$ les $\OO$-modules inversibles associ\'es. Soit $g : Y'\rr Y$ un morphisme de sch\'emas sur $X$. D\'esignons par }
$$
\psi : f_{\star}\oo_{Y} \rr f'_{\star}\oo_{Y'}
$$
{\it le morphisme de $\OO$-alg\`ebres associ\'e \`a $g$. Alors

\n i) Si $\psi$ est injectif, on a $\psi(N) \subset N'$.

\n ii) Si $\psi$  n'est pas injectif, et si $X$ est normal int\`egre, alors $g = s f' $ o\`u $s$ est une section de $f$.}
\medskip

\' Ecrit par blocs, le morphisme $\psi : \OO \oplus N \rr \OO \oplus N'$ prend la forme suivante
$$
\psi = \begin{pmatrix}
1 & \psi_{1}\\
0 & \psi_{0}
\end{pmatrix}
$$
{\it i)}\; Il s'agit de voir que la forme lin\'eaire $\psi_{1} : N \rr \OO$ est nulle. Notons d'abord que l'application \;  $\psi_{0}: N \rr N'$ est injective, car si un \'el\'ement $\alpha \in N$ est tel que $\psi_{0}(\alpha) = 0$, alors on a \;  $-\psi_{1}(\alpha) + \alpha \in {\rm Ker}(\psi)$, donc $\alpha = 0$.

\n Soit $\alpha$ une section de $N$, de sorte que $\alpha^2$ est une section de $\OO$ et que, par suite $\psi(\alpha^2) = \alpha^2$. Comme $\psi$ respecte le produit, on a
$$
\alpha^2 = \psi(\alpha^2) = \psi(\alpha)^2 = (\psi_{1}(\alpha)+ \psi_{0}(\alpha))^2 = \psi_{1}(\alpha)^2 + \psi_{0}(\alpha)^2 \; + \; 2\psi_{1}(\alpha)\psi_{0}(\alpha).
$$
Le dernier terme, $2\psi_{1}(\alpha)\psi_{0}(\alpha)$,\, est la composante dans $N'$ ; il est donc nul. Comme $\psi_{0}$ est injectif, on voit que l'on a, pour tout $\alpha \in N$,  $\psi_{1}(\alpha)\alpha = 0$. Mais $N$ \'etant un module inversible, cette relation entra\^ine la nullit\'e de $\psi_{1}$.
\medskip

{\it ii)}\; Supposons que $X$ soit normal int\`egre, et que $\psi$ ne soit pas injectif ; son noyau est donc g\'en\'eriquement de rang 1; par suite, la $\OO$-alg\`ebre ${\mathcal B} = {\rm Im}(\psi) \subset f'_{\star}\oo_{Y'}$ est finie, g\'en\'eriquement de rang 1, et elle est sans torsion puisque contenue dans  $f'_{\star}\oo_{Y'}$ ; comme $\OO$ est normale, ${\mathcal B}  = \OO$. D'o\`u le r\'esultat.
\bigskip

Voir {\bf 2.2.} pour une description plus compl\`ete de $g$ dans le cas {\it i)}.
\bigskip

\noindent {\bf 1.7.}\; Rev\^etements localement isomorphes
\bigskip

\noindent {\bf Lemme}\; {\it Soient $f_{1} : Y_{1} \rightarrow X$ et $f_{2} : Y_{2} \rightarrow X$ deux rev\^etements doubles de $X$, et $N_{1}$ et $N_{2}$ les $\OO$-modules inversibles associ\'es. On suppose que le lieu de diramation de $f_{1}$ et celui de  $f_{2}$ sont des diviseurs. Alors, les propri\'et\'es suivantes sont \'equivalentes}

\noindent {\it i) Il existe un morphisme fid\`element plat quasi-compact $p : X' \rightarrow X$ et un isomorphisme}
$$ X'\times_{X}Y_{1} \; \widetilde{\longrightarrow}\; X'\times_{X}Y_{2} \leqno{(1.7.1)}$$

\noindent {\it ii) Les lieux de diramation de $f_{1}$  et de $f_{2}$  sont \'egaux.}

\medskip

${\it i)} \Rightarrow {\it ii)}$. Comme la trace commute \`a un isomorphisme, la donn\'ee de (1.7.1)  \'equivaut \`a la donn\'ee d'un isomorphisme de ${\mathcal O}_{X'}$-modules inversibles
$$
p^{\star}(N_{1})\; \stackrel{\omega}{\longrightarrow}\; p^{\star}(N_{2})
$$
compatible avec les structures multiplicatives   $\mu_{i} : N_{i}^{\otimes 2} \longrightarrow \OO$ ; comme ces applications sont suppos\'ees  injectives, l'isomorphisme
$$
p^{\star}(N_{1}^{\otimes 2})\; \stackrel{\omega^{\otimes 2}}{\longrightarrow}\; p^{\star}(N_{2}^{\otimes 2})
$$
se descend de $X'$ \`a $X$ ; par suite, $\mu_{1}$  et $\mu_{2}$ ont m\^eme image.

${\it ii)} \Rightarrow {\it i)}$. Par hypoth\`ese, il existe un isomorphisme $\theta : N_{1}^{\otimes 2}\quad  \widetilde{\longrightarrow}\quad  N_{2}^{\otimes 2}$ tel que $\mu_{2}\circ \theta = \mu_{1}$. 

\noindent Consid\'erons le foncteur $F$ d\'efini, pour tout $X$-sch\'ema $q : Z \rightarrow X$, par 
$$
F(Z) = \{ \omega \in {\rm Hom}_{{\mathcal O}_{Z}}(q^{\star}(N_{1}), q^{\star}(N_{2}))\; {\rm tel \, que}\; \omega^{\otimes 2} = q^{\star}(\theta)\}
$$
Ce foncteur est  repr\'esentable par un rev\^etement \'etale de rang deux. En effet, posons $L = {\mathcal Hom}(N_{1}, N_{2})$ ; c'est un module inversible dont le carr\'e est muni de l'isomorphisme $ \lambda : L^{\otimes 2}  \quad \widetilde{\longrightarrow}\quad  \OO$, donn\'e par $\theta$.

Alors le rev\^etement  double\;  $p : X' = {\mathcal Spec}(\OO \oplus L) \; \rightarrow X$, associ\'e \`a $\lambda$, repr\'esente le foncteur en question ({\bf 1.5}) ; par suite on a isomorphisme $ X'\times_{X}Y_{1} \quad \widetilde{\longrightarrow}\quad X'\times_{X}Y_{2}.$
\vspace{1cm}

%%%%%%%%%%%%%%%%%%%%%%%%

\noindent {\bf 2.\; Exemples}
\bigskip

\noindent {\bf 2.1.}\; Construction par pincement
\bigskip

 Soit $f' : Y' \rightarrow X$ un rev\^etement double, et $D \subset X$ un diviseur (de Cartier) effectif sur $X$. Posons $D' = f'^{-1}(D)$ ; c'est un diviseur effectif sur $Y'$, et le morphisme  $\bar{g} : D' \rightarrow D$ induit par $f'$, est un rev\^etement double. En \og pin\c{c}ant \fg\; $Y'$ le long de $\bar{g}$, on obtient  un sch\'ema $Y$ et un diagramme commutatif
  $$
  \xymatrix{D' \ar[r]^{\bar{g}} \ar[d] & D \ar[d] \ar[dr] &\\
 Y' \ar[r]^{g} \ar@/_1pc/[rr]_{f'}&Y \ar[r]^{f} &X }
 $$
 o\`u le carr\'e est cocart\'esien, o\`u $ f' = f\circ g$, et o\`u les morphismes verticaux sont des immersions ferm\'ees.Voir [Fer. 2], thm. 7.1.B.(L'hypoth\`ese {\it iii)} de {\it loc.cit.} se r\'eduit \`a ceci : pour tout $x \in D \subset X$, la fibre $f'^{-1}(x)$ est contenue dans un ouvert affine de $Y'$ ; or, c'est \'evident ici puisque le morphisme $f'$ est affine). 
 \medskip
 
 Montrons que $f : Y \rightarrow X$ est un rev\^etement double.
 Montrons d'abord que $f$ est un morphisme affine : soit $U$ un ouvert affine de $X$ ; alors $f^{-1}(U)$ est la somme des sch\'emas $D \cap U$ et $f'^{-1}(U)$, amalgam\'ee le long de $f'^{-1}(D \cap X)$ :
 $$
 \begin{CD}
 f'^{-1}(D \cap X) @>>> D\cap X\\
 @VVV @VVV\\
 f'^{-1}(U) @>>> f^{-1}(U)
 \end{CD}
 $$
 
  Or, ces trois derniers sch\'emas sont affines ; donc $f^{-1}(U)$ est affine (loc.cit. Thm 5.1). Il reste \`a v\'erifier que $f$ est localement libre de rang deux. Consid\'erons les images directes sur $X$ des faisceaux d'anneaux des diff\'erents sch\'emas qui interviennent ; on peut omettre  les symboles d'image directe $f_{\star}$, etc. puisque les morphismes sont affines. La propri\'et\'e \`a v\'erifier \'etant locale, on peut supposer que tous les sch\'emas sont affines.
  
  \noindent Notons $\mathcal{L}$ et $\mathcal{M}$ les $\OO$-modules quasi-coh\'erents conoyaux des morphismes $ \mathcal{O}_{D}\longrightarrow \mathcal{O}_{D'}$, et $\mathcal{O}_{Y}\longrightarrow \mathcal{O}_{Y'}$. Par d\'efinition d'un carr\'e cocart\'esien, dans le diagramme suivant, les lignes sont exactes et l'application $\mathcal{M} \rightarrow \mathcal{L}$ est un isomorphisme.
  $$
  \xymatrix{0   & \mathcal{L} \ar[l] &\mathcal{O}_{D'} \ar[l] & \mathcal{O}_{D} \ar[l] &0 \ar[l]\\
  0    & \mathcal{M} \ar[u] \ar[l] &\mathcal{O}_{Y'} \ar[u] \ar[l] & \mathcal{O}_{Y} \ar[u] \ar[l] &\ar[l] 0}
  $$
 Comme $D' \rightarrow D$ est un rev\^etement de rang deux,  $ \mathcal{L} $ est un $ \mathcal{O}_{D}$-module inversible ; mais $D$ est un diviseur de $X$, donc
 $$
 {\rm dim.proj}_{ \mathcal{O}_{X}}( \mathcal{L}) = 1.
 $$
 (Les notions de module projectif et de dimension projective ont bien un sens puisque $X$ est affine). Puisque $ \mathcal{O}_{Y'}$ est localement libre de rang 2 sur $\OO$, on d\'eduit de l'isomorphisme $\mathcal{M} \tilde{\rightarrow} \mathcal{L}$, et de  l'exactitude de la suite inf\'erieure  du diagramme, que $ \mathcal{O}_{Y}$ est projectif sur  
 $\OO$. D'autre part,  c'est un $\OO$-module de type fini puisque $\mathcal{M}$ est de pr\'esentation finie et que $\mathcal{O}_{Y'}$  est de type fini ([Bour] AC I \S2.8).
 Finalement, comme $f' : Y' \rightarrow Y$ est un isomorphisme au dessus de l'ouvert $X - D$, on peut bien conclure que le rang du  $\OO$-module projectif $\oo_{Y}$ est \'egal \`a 2. \hspace{10cm}$\Box$
 \medskip
 
 Voici une description alg\'ebrique de  $\A = f_{\star}({\mathcal O}_{Y})$ : posons $\A' = f'_{\star}({\mathcal O}_{Y'}) = \OO \oplus N'$ ; notons ${\mathcal J} \subset \OO$ l'id\'eal (inversible) de $D$ dans $X$ ; alors $\A  = \OO \oplus {\mathcal J}N'  \subset \A'$.
 \bigskip
 
 \n En fait, on va voir que \emph{tout} morphisme \og raisonnable \fg \, entre rev\^etements doubles est un pincement au sens pr\'ec\'edent.
 \bigskip
 
 \n {\bf 2.2. Proposition}\; {\it Soient $f : Y \rr X$  et $f' : Y' \rr X$ deux rev\^etements doubles de $X$, et $(N, \mu)$, $(N', \mu')$ les modules inversibles, et les multiplications associ\'es. On consid\`ere un $X$-morphisme $g : Y' \rr Y$. On suppose que le morphisme de $\OO$-alg\`ebres associ\'e \`a $g$}
 $$
 \psi : \A=\OO \oplus N \; \rr \; \OO \oplus N'=\A'
 $$
 {\it est \emph{injectif}, de sorte que $\psi(N) \subset N'$}\, ({\bf 1.6}).
 
 \n{\it  Soit $D \subset X$ le diviseur de Cartier d\'efini par la suite exacte d\'eduite de $\psi$}
 $$
 0 \rr {\mathcal Hom}(N', N) \rr \OO \rr \oo_{D} \rr 0
 $$
 {\it Alors, $Y$ s'obtient par pincement de $Y'$ le long du morphisme  $f'^{-1}D = D' \, \rr \, D$.}
 \medskip
 
 Identifions $N$ \`a son image $\psi(N) \subset N'$, ainsi que $\A$ \`a son image dans $\A'$. Notons d'abord que l'id\'eal de $D$ dans $X$ est \'egal \`a $\mathcal{J} = {\rm Ann}_{\OO}(N'/N)$. Puisque $N$ et $N'$ sont des $\OO$-modules inversibles, on a m\^eme l'\'egalit\'e $\mathcal{J}N' = N$, et, par suite, $\mathcal{J}\A' = \mathcal{J} \oplus \mathcal{J}N' = \mathcal{J} \oplus N$ . De plus, on v\'erifie imm\'ediatement  que le conducteur de $\psi$, c'est-\`a-dire l'id\'eal ${\rm Ann}_{\A}(\A' / \A)$, est  \'egal \`a $\mathcal{J}\oplus N$, et donc aussi \`a $\mathcal{J}\A'$.
 
 \n Par d\'efinition du conducteur, le carr\'e 
 $$
 \begin{CD}
 \A' / \mathcal{J}\A' = \A' / \mathcal{J}\oplus N @<<< \A / \mathcal{J}\oplus N \simeq \OO / \mathcal{J} \\
 @AAA @AAA\\
 \A' @<<< \A
 \end{CD}
 $$
 est cocart\'esien ({\it i.e.} il fait de $\A$ le produit fibr\'e \'evident), ce qui est une autre fa\c{c}on de dire qu'en passant aux sch\'emas, le carr\'e obtenu
$$
\begin{CD}
D' @>>> D\\
@VVV @VVV \\
Y' @>>g> Y
\end{CD}
$$
est un pincement. 
\vspace{1cm}

\n 

\noindent {\bf 2.3.}\; Rev\^etement standard
\bigskip

\noindent {\bf 2.3.1.\; D\'efinition} \; {\it Soit $Z$ un diviseur de Cartier effectif sur un sch\'ema $X$. On appelle rev\^etement standard de $X$ associ\'e \`a $Z$ le sch\'ema obtenu en recollant deux copies de $X$ le long de $Z$.}
\bigskip

Ce type de recollement est un cas particulier du pr\'ec\'edent, o\`u $Y' = X \sqcup X$, et son existence se trouve d\'ej\`a dans [Anan.] Prop.1.1.1.
\bigskip

Soit $f : Y \rightarrow X$ le rev\^etement standard associ\'e au diviseur $Z \subset X$ ; par d\'efinition, on a un carr\'e de somme amalgam\'ee (deux copies de $X$ amalgam\'ees en $Z$)
$$
\begin{CD}
Z @>>> X\\
@VVV  @VVV\\
X @>>> Y
\end{CD}\leqno{(2.3.2)}
$$
Posant, comme plus haut, $\A = f_{\star}(\mathcal{O}_{Y})$, on obtient une suite exacte de $\OO$-modules
$$
0\; \longrightarrow\; \A \; \stackrel{\iota}{\longrightarrow}\;  \OO \times \OO \; \xrightarrow{(x,y) \mapsto {\rm cl}(y-x)}\;  \mathcal{O}_{Z} \; \longrightarrow\;  0 .
$$
Comme $\mathcal{O}_{Z} = \OO / \mathcal{J}$, o\`u $\mathcal{J}$ est un id\'eal inversible, on voit que $\A$ est localement libre, et de rang 2. L'application $\iota$ est un isomorphisme au dessus de l'ouvert sch\'ematiquement dense $X - Z$ ; par suite, la trace se calcule comme pour les \'el\'ements de $\OO \times \OO $ : c'est la somme des composantes ; le noyau de la trace  $N = {\rm Ker}(\A \stackrel{{\rm Tr}}{\rightarrow} \OO)$ est donc isomorphe \`a $\mathcal{J}$ :
$$
\mathcal{J} \quad  \widetilde{\longrightarrow} \quad N, \hspace{1cm} x \mapsto (x, -x).
$$
Le diviseur de diramation est d\'efini par l'id\'eal $\mathcal{I}$ image de $N^{\otimes 2} \rightarrow \OO$, soit
$$
\mathcal{I}\; = \; \mathcal{J}^2.
$$
\medskip

\noindent {\bf 2.4.} \; On peut d\'ecrire $\A$ plus simplement : on munit $\OO \oplus \mathcal{J}$ de la structure de $\OO$-alg\`ebre pour laquelle le produit de deux \'el\'ements du second facteur $\mathcal{J}$ est ce produit vu comme \'el\'ement de $\OO$ ; alors, on a un isomorphisme de $\OO$-alg\`ebres
$$
\OO \oplus \mathcal{J} \simeq \A\;  \subset \;  \OO \times \OO, \qquad  t \oplus u  \mapsto (t-u, t+u).
$$
La d\'emarche est r\'eversible ; cela montre que pour tout id\'eal inversible $\mathcal{J} \subset \OO$, le spectre de l'alg\`ebre $\OO \oplus \mathcal{J}$ est le rev\^etement standard associ\'e au diviseur ${\mathcal Spec}(\OO / \mathcal{J})$.

\noindent Lorsque $X = {\rm Spec}(R)$ et que $\mathcal{J} = uR$, alors $Y$ est isomorphe \`a ${\rm Spec}(R[T]/(T^2-u^2))$.
\bigskip

\noindent {\bf 2.4.1. \, Proposition}\; {\it Soit $f : Y \rightarrow X$ un rev\^etement de rang deux dont le lieu de diramation est un diviseur. Alors $f$ est standard si et seulement si il admet une section.}
\bigskip

L'existence d'une section est clairement n\'ecessaire ; montrons qu'elle est suffisante. Soit  $\mu : N^{\otimes 2} \rightarrow \OO$ la multiplication associ\'ee \`a $Y$. Soit $s : X \rightarrow Y$  une section de $f$. Consid\'erons, comme en {\bf 1.5}, l'application $f^{\star}(N) \rightarrow \mathcal{O}_{Y}$ adjointe de l'inclusion ; son image r\'eciproque par la section $s$ donne une application $\OO$-lin\'eaire $\omega : N \rightarrow \OO$, dont le carr\'e est \'egal \`a $\mu$ ; la remarque ci-dessus montre que $Y$ est standard.$\Box$
\medskip

Cela montre que pour tout rev\^etement de rang deux, $f : Y \rightarrow X$, \`a diramation port\'ee par un diviseur,  le rev\^etement, d\'eduit par changement de base, $Y \times_{X} Y \longrightarrow Y$ est standard : il est obtenu par recollement de deux copies de $Y$ le long du ferm\'e d'id\'eal $\mathcal{I} \oplus N \subset  f_{\star}(\mathcal{O}_{Y})$, c'est-\`a-dire le long du lieu de diramation $\Delta$ vu comme ferm\'e de $Y$ (cf. {\bf 1.3}).
\vspace{1cm}

%%%%%%%%%%%%%%%%%%%%%%%%%%
\noindent {\bf 3.}\; {\bf Le rev\^etement  double }$({\bf P}_{1})^n / \mathfrak{A}_{n}\; \longrightarrow \; {\bf P}_{n}$
\bigskip

\n Une des formes du \emph{ th\'eor\`eme des polyn\^omes  sym\'etriques \'el\'ementaires} ([Bour], A IV.58, et TG VIII.22) conduit \`a un isomorphisme $$({\bf P}_{1})^n / \mathfrak{S}_{n}\; \; \tilde{\longrightarrow} \; {\bf P}_{n}$$
On s'int\'eresse ici au passage au quotient par le groupe altern\'e $\A_{n}$, qui est d'indice deux dans $\mathfrak{S}_{n}$, et qui donne donc lieu \`a un rev\^etement de degr\'e 2 de ${\bf P}_{n}$. Dans ce \S, on va d\'eterminer  le faisceau inversible associ\'e \`a ce rev\^etement, au sens de la d\'efinition {\bf 1.2}, ainsi que sa multiplication. L'analogue affine est connu et fa\c{c}ile (voir {\bf 3.7.2}) ; mais le r\'esultat global vis\'e ici requiert des v\'erifications parfois fastidieuses que, cependant, nous avons voulu ne pas escamoter.\\

Dans ce paragraphe, tous les sch\'emas sont sur Spec($K$), o\`u $K$ est un corps de caract\'eristique $\neq 2$, et ce sch\'ema de base est sous-entendu.\\

\n {\bf 3.1.}\; Soit $V$ un espace vectoriel de dimension finie sur le corps $K$, et ${\bf P}(V)$ l'espace projectif associ\'e ; on note $\alpha : V \rightarrow L$ le quotient inversible fondamental. Ici, et plus bas, lorsqu'il n'y a pas d'ambigu\"{\i}t\'e, on \'ecrit simplement $V$ \`a la place de son image r\'eciproque $V_{\pp(V)}$, etc. 
Soit $\pp(V)^n$ le produit de $n$ copies de $\pp(V)$, et $p_{i} : \pp(V)^n \rightarrow \pp(V)$ la projection sur le facteur d'indice $i$ ; sur le sch\'ema produit on obtient, par image r\'eciproque, des quotients inversibles
$$
\alpha_{i}= p_{i}^{\star}(\alpha) : V \longrightarrow L_{i} = p_{i}^{\star}(L).
$$
On note ${\sf TS}^n(V)$ l'espace des tenseurs sym\'etriques. 

\n L'application compos\'ee 
$$
{\sf TS}^n(V) \; \subset \; V^{\otimes n} \; \xrightarrow{\alpha_{1}\otimes \cdots \otimes \alpha_{n}} L_{1}\otimes \cdots \otimes L_{n}
$$
est surjective (il s'agit d'applications lin\'eaires entre modules sur $\pp(V)^n$). Pour le voir, on peut supposer effectu\'e  un changement de base de la forme ${\rm Spec}(K') \rightarrow \pp(V)^n$, o\`u $K'$ est une extension de $K$ ; les $L_{i}$ sont alors des vectoriels de rang 1 ; si $\alpha_{1}(x_{1}) \otimes\ldots \otimes \alpha_{n}(x_{n})$ est un \'el\'ement non nul de $L_{1}\otimes \cdots \otimes L_{n}$, il existe  des \'el\'ements inversibles $\lambda_{i} \in K'$ tel que, dans $L_{i} \simeq K'$, on ait  $\alpha_{i}(x_{i}) = \lambda_{i} \alpha_{i}(x_{1})$ ; par suite,  on a 
$$
\alpha_{1}(x_{1}) \otimes\ldots \otimes \alpha_{n}(x_{n}) = \lambda_{1}\lambda_{2}\ldots \lambda_{n}\alpha_{1}(x_{1}) \otimes\ldots \otimes \alpha_{n}(x_{1}).
$$ 
C'est l'image de l'\'el\'ement sym\'etrique $ \lambda_{1}\lambda_{2}\ldots \lambda_{n}x_{1}^{\otimes n}$. 
\bigskip

\n {\bf 3.2.}\; On suppose maintenant que  $V$ est de rang 2, donc que ${\bf P}(V)$ est une droite projective ; le sch\'ema ${\bf P}({\sf TS}^n(V))$ est isomorphe \`a ${\bf P}_{n}$ puisque ${\sf TS}^n(V)$ est de rang $n+1$ ; pour faire court on \'ecrira parfois ${\bf P}_{n}$ \`a la place de ${\bf P}({\sf TS}^n(V))$, m\^eme sans choix de base explicit\'e. La surjectivit\'e montr\'ee ci-dessus donne donc lieu \`a un morphisme
$$
\pi :  {\bf P}(V)^n \longrightarrow \; {\bf P}({\sf TS}^n(V)) \simeq {\bf P}_{n} .
$$
Il est caract\'eris\'e par l'existence d'un isomorphisme
$$
\pi^{\star}({\mathcal O}_{{\bf P}_{n}}(1)) \; \simeq \; L_{1}\otimes \cdots \otimes L_{n} \leqno{(3.2.1)}
$$
rendant commutatif le diagramme suivant d'applications lin\'eaires entre  modules sur ${\bf P}(V)^n $ :
$$
\xymatrix{
{\sf TS}^n(V) \ar@{=}[d]  \ar[rr]&& \pi^{\star}({\mathcal O}_{{\bf P}_{n}}(1))\ar[d]^{\wr} \\
{\sf TS}^n(V) \ar[r] & V^{\otimes n} \ar[r] &L_{1}\otimes \cdots \otimes L_{n}
 } \leqno{(3.2.2)}
$$
\medskip

\n Rappelons la d\'efinition  du morphisme $\pi$ en termes de coordonn\'ees homog\`enes.

\n Soit $\{e_{0}, e_{1}\}$ une base de $V$. Pour toute partie $I \subset \{1, \ldots , n\}$, on note $e_{I} \in V^{\otimes n}$ le produit
$$
e_{I} = v_{1}\otimes \cdots \otimes v_{n}, \quad \textrm{o\`u}\quad  v_{i} = \left \{
\begin{array}{cc}
e_{0} & \textrm{si} \; i\in I\\
e_{1} & \textrm{sinon}
\end{array}
\right.
$$
Pour $p = 0, 1, \ldots , n$, on pose
$$
{\bf e}_{p} \; = \; \sum_{| I | = p} e_{I}.
$$
C'est un tenseur sym\'etrique ; mieux : $\{ {\bf e}_{0}, \ldots , {\bf e}_{n} \}$ est une base du $K$-espace vectoriel ${\sf TS}^n(V)$ ( [Bour] A IV 5.5 Prop. 4). Enfin, on a,  en notant $T$ une ind\'etermin\'ee,
$$
(e_{0} + Te_{1})^{\otimes n} \; = \; {\bf e}_{n} + {\bf e}_{n-1}T + \cdots + {\bf e}_{0}T^n.\leqno{(3.2.3)}
$$
Consid\'erons $n$ points de ${\bf P}_{1}(K)$ d\'etermin\'es par les formes lin\'eaires surjectives $\alpha_{i} : V \rightarrow K$ ; on \'ecrit  $x_{i} = \alpha_{i}(e_{0})$, et  $y_{i} = \alpha_{i}(e_{1})$ ; posons alors, pour une partie $I \subset \{1, \ldots , n\}$, 
$$
z_{I} := {\alpha_{1}\otimes \cdots \otimes \alpha_{n}}(e_{I}) = \prod_{i \in I}x_{i} \prod_{i \notin I} y_{i},
$$
et, pour $p = 0, 1, \ldots , n$, 
$$
{\bf z}_{p} \; = \;  \sum_{| I | = p} z_{I}
$$
L'\'egalit\'e (3.2.3) conduit donc \`a l'\'egalit\'e suivante entre polyn\^omes dans $K[T]$
$$
(x_{1} + y_{1}T)\cdots (x_{n} + y_{n}T) \quad = \quad {\bf z}_{n} + {\bf z}_{n-1}T + \cdots + {\bf z}_{0}T^n.
$$
Le morphisme $\pi :  {\bf P}_{1}^n \longrightarrow \; {\bf P}_{n}$ s'\'ecrit alors
$$
(x_{1} : y_{1}) , \ldots , (x_{n} : y_{n}) \; \longmapsto \; ({\bf z}_{0} :  {\bf z}_{1} : \cdots : {\bf z}_{n} ).
$$
\bigskip

\n {\bf 3.3.}\; On va d\'efinir un module inversible ${\sf M}$ sur ${\bf P}(V)^n$, et montrer qu'il est isomorphe \`a l'image r\'eciproque $\pi^{\star}{\sf M}_{0}$ d'un module inversible ${\sf M}_{0}$ d\'efini sur $\pp^n$ ; on montrera ensuite que le module inversible associ\'e au rev\^etement double consid\'er\'e est le dual de ${\sf M}_{0}$. 
\medskip

\n D\'esignons par ${\mathcal C} = {\mathcal C}^2_{n}$ l'ensemble des parties \`a 2 \'el\'ements de $\{1, \ldots , n\}$. Pour $J \in {\mathcal C}$, si on \'ecrit  $J = \{i, j\}$, avec $i < j$, on pose
$$
{\sf M}_{J} = {\mathcal Hom}(\W^2 V, L_{i}\otimes L_{j})
$$
C'est un module inversible sur $\pp(V)^n$. 

\n On pose
$$
{\sf M}\; = \bigotimes_{J \in {\mathcal C}} {\sf M}_{J}.
$$
(Disons qu'on a choisi l'ordre lexicographique sur ${\mathcal C}$). 

\n Montrons que le module inversible ${\sf M}$ sur ${\bf P}(V)^n$ provient, via le morphisme
$$\pi :  {\bf P}(V)^n \longrightarrow \; {\bf P}({\sf TS}^n(V))$$ d'un module inversible ${\sf M}_{0}$  sur la base ${\bf P}({\sf TS}^n(V))$, laquelle va \^etre not\'ee simplement ${\bf P}_{n}$. 

\n La d\'efinition de $\pi$ est li\'ee \`a l'isomorphisme (3.2.1)
$$
\pi^{\star}({\mathcal O}_{{\bf P}_{n}}(1)) \; \simeq \; L_{1}\otimes \cdots \otimes L_{n}.
$$
Or, parmi les  $\frac{n(n-1)}{2}$ parties \`a deux \'el\'ements de $\{1, \ldots , n\}$ il y en a $n-1$  contenant un \'el\'ement fix\'e $i$ ; en utilisant l'isomorphisme canonique ${\sf M}_{J} = {\mathcal Hom}(\W^2 V, L_{i}\otimes L_{j}) \simeq (\W^2V)^{\otimes - 1} \otimes L_{i}\otimes L_{j}$, on trouve donc un isomorphisme
$$
{\sf M}\; = \bigotimes_{J \in {\mathcal C}} {\sf M}_{J} \quad  \widetilde{\rr} \quad  \bigotimes_{J \in {\mathcal C}}(\W^2V)^{\otimes - 1} \otimes L_{i}\otimes L_{j}  \quad  \widetilde{\rr} \quad (\W^2V)^{\otimes - \frac{n(n-1)}{2}} \otimes (L_{1}\otimes \cdots \otimes L_{n})^{\otimes (n-1)} .
$$

Compte-tenu de (3.2.1), on obtient, finalement, un isomorphisme 
$$
{\sf M} \quad \widetilde{\rr} \quad  (\W^2V)^{\otimes - \frac{n(n-1)}{2}} \otimes \pi^{\star}({\mathcal O}_{{\bf P}_{n}}(n-1)).
$$
Posons donc 
$$\boxed{{\sf M}_{0} = (\W^2V)^{\otimes - \frac{n(n-1)}{2}} \otimes {\mathcal O}_{{\bf P}_{n}}(n-1)} \leqno{(3.3.1)}
$$
de sorte qu'on a d\'egag\'e un isomorphisme
$$
\eta : {\sf M} \rr \pi^{\star}({\sf M}_{0}).\leqno{(3.3.2)}
$$
\bigskip

\n {\bf 3.4.}\; Les permutations et leur signature joueront \'evidemment un r\^ole central dans la suite ; il faut donc d'abord pr\'eciser l'op\'eration (\`a gauche) de $\mathfrak{S}_{n}$ sur diff\'erents modules d\'efinis sur $\pp(V)^n$.

\n L'op\'eration de $\mathfrak{S}_{n}$ sur le sch\'ema $\pp(V)^n$ lui-m\^eme est caract\'eris\'ee, en termes des projections $p_{i}$, par la relation suivante : pour tout $\sigma \in \mathfrak{S}_{n}$, on a
$$
\boxed{p_{i}\circ \sigma = p_{\sigma^{-1} i}}
$$
On trouve donc
$$
\sigma^\star(V \stackrel{\alpha_{i}}{\longrightarrow} L_{i}) = \sigma^\star p_{i}^\star(V \stackrel{\alpha}{\longrightarrow} L) = p_{\sigma^{-1}i}^\star(V \stackrel{\alpha}{\longrightarrow} L) = (V \stackrel{\alpha_{\sigma^{-1}i}}{\longrightarrow} L_{\sigma^{-1}i}) 
$$
En particulier, on a
$$
\boxed{\sigma^\star(L_{i}) = L_{\sigma^{-1}i} .}
$$
\medskip

\n Fixons une permutation $\sigma$ ; on a un isomorphisme
$$
\theta_{J} : \sigma^\star {\sf M}_{J} \quad \widetilde{\rr}\quad {\sf M}_{\sigma^{-1}J}\, .
$$
Il fait intervenir l'isomorphisme de commutativit\'e $L_{\sigma^{-1}i}\otimes L_{\sigma^{-1}j} \simeq L_{\sigma^{-1}j}\otimes L_{\sigma^{-1}i}$ lorsque $\sigma^{-1}i > \sigma^{-1}j$.

On d\'esigne par
$$
\theta : \sigma^\star {\sf M} \quad \widetilde{\rr}\quad {\sf M}\leqno{(3.4.1)}
$$
l'isomorphisme obtenu en composant le produit tensoriel des $\theta_{J}$ avec l'isomorphisme canonique de permutation des facteurs \, $ \bigotimes_{J \in {\mathcal C}} {\sf M}_{\sigma^{-1}J} \;  \simeq \; \bigotimes_{J \in {\mathcal C}} {\sf M}_{J}$ [Bour] A II.104.

\n Concr\`etement, l'application $\theta$ est $\sigma$-lin\'eaire, {\it i.e}  c'est une application additive ${\sf M} \rr {\sf M}$, telle que, pour des sections locales  $s$ et $m$, on ait $\theta(sm) = \sigma (s)\theta(m)$.
\medskip

Explicitons la compatibilit\'e de $\eta$ avec la permutation $\sigma$ : le carr\'e suivant est commutatif
$$
\xymatrix{
 \sigma^\star{\sf M}\ar[d]_{\theta} \ar[r]^{\sigma^\star(\eta)} & \sigma^{\star}\pi^{\star}({\sf M}_{0}) \ar@{=}[r] &\pi^{\star}({\sf M}_{0}) \ar[d]^{\sigma} \\
{\sf M} \ar[rr]_{\eta}&& \pi^{\star}({\sf M}_{0})
 }
$$

En effet, dans le carr\'e (3.2.2), les applications horizontales (de source le module $\mathfrak{S}_{n}$-invariant ${\sf TS}^n(V)$) sont surjectives ; par suite l'isomorphisme (3.2.1) est compatible aux permutations, et la commutativit\'e du carr\'e ci-dessus en d\'ecoule.
\bigskip

\n {\bf 3.5.}\; On va maintenant d\'efinir une section  $w  :  \oo_{\pp(V)^n} \; \rr \; {\sf M}$, et v\'erifier qu'elle est invariante sous le groupe altern\'e, et que son carr\'e est invariante sous $\mathfrak{S}_{n}$.
\medskip

\n Pour tout $J \in {\mathcal C}$, consid\'erons la  section
$$
w_{J} : \oo_{\pp(V)^n} \; \rr \; {\sf M}_{J}, 
$$
associ\'ee \`a l'application 
$$
\W^2V \rr L_{i}\otimes L_{j}, \qquad  x\wedge y \longmapsto \alpha_{i}(x)\otimes \alpha_{j}(y) -\alpha_{i}(y)\otimes \alpha_{j}(x).\leqno{(3.5.1)}
$$
Il faut remarquer que le carr\'e
$$
\xymatrix{
\oo_{\pp(V)^n} \ar[d]_{\sigma}  \ar[r]^{\sigma^\star(w_{J})}& \sigma^\star{\sf M}_{J}\ar[d]^{\theta_{J}} \\
\oo_{\pp(V)^n} \ar[r]_{w_{\sigma^{-1}J}} & \; {\sf M}_{\sigma^{-1}J}
 }
$$
est commutatif ou anti-commutatif selon que l'on a  $\sigma^{-1} i < \sigma^{-1} j$, ou bien $\sigma^{-1}i > \sigma^{-1} j$.

Le support du diviseur div$(w_{J})$ est le \emph{ferm\'e d'\'egalit\'e}  de $\alpha_{i}$ et de $\alpha_{J}$, ou, si l'on pr\'ef\`ere, le ferm\'e form\'e des points dont les coordonn\'ees d'indices $i$ et $j$ sont \'egales ; ce sous-sch\'ema est invariant par la transposition $(i\, j)$, mais elle transforme la section $w_{J}$ en son oppos\'ee $- w_{J}$ ; pour d\'egager le diviseur cherch\'e on ne peut donc malheureusement  pas s'en tenir \`a l'intuition purement g\'eom\'etrique et identifier un diviseur \`a son support .
\medskip

\n Par produit tensoriel des $w_{J}$, on obtient donc une application
$$
w  :  \oo_{\pp(V)^n} \; \rr \; {\sf M}
$$
telle que le carr\'e suivant soit commutatif

$$
\xymatrix{
\oo_{\pp(V)^n} \ar[d] _{\sigma} \ar[r]^{\sigma^\star(w)}& \sigma^\star{\sf M}\ar[d]^{\theta} \\
\oo_{\pp(V)^n} \ar[r]_{\varepsilon(\sigma)w} & {\sf M}
 }\leqno{(3.5.2)}
$$
o\`u  on a not\'e $\varepsilon(\sigma)$ la signature de la permutation $\sigma$. 
\bigskip

\n {\bf 3.6.}\; Finalement, on a d\'efini une section, sur le sch\'ema produit ${\bf P}(V)^n$,
$$
\omega = \eta \circ w : \mathcal{O}_{\pp(V)^n} \; \rr \; \pi^{\star}({\sf M}_{0})
$$
qui est invariante par permutation paire, au sens o\`u le carr\'e suivant est commutatif

$$
\xymatrix{
\oo_{\pp(V)^n} \ar[d] _{\sigma} \ar[r]^{\sigma^\star(\omega)}& \pi^\star{\sf M}_{0}\ar[d]^{\sigma} \\
\oo_{\pp(V)^n} \ar[r]_{\varepsilon(\sigma)\omega} & \pi^{\star}{\sf M}_{0}
 }\leqno{(3.6.1)}
$$
\medskip

\n All\'egeons maintenant en $f : Y \rr X$ la notation du morphisme  qui est en cause :
$$
Y = \pp(V)^n / \mathfrak{A}_{n} \quad \stackrel{f}{\rr}\quad X =  \pp(V)^n / \mathfrak{S}_{n} = \pp({\sf TS}^n(V))
$$
L'invariance de la section $\omega$ par permutation paire entra\^ine qu'elle provient d'une section d\'efinie sur $Y$, et que son carr\'e $\omega^{\otimes 2}$ se descend \`a $X$ ; il existe donc une application $\oo_{Y}$-lin\'eaire, not\'ee par la m\^eme lettre,
$$
\omega : \oo_{Y} \; \rr \; f^{\star}({\sf M}_{0})\leqno{(3.6.2)}
$$
telle que $\omega^{\otimes 2} : \oo_{Y} \; \rr \; f^{\star}({\sf M}^{\otimes 2}_{0})$ provienne d'une application d\'efinie sur $X$.

\n Posons alors
$$
\boxed{N \; = \; {\mathcal Hom}_{\oo_{X}}({\sf M}_{0}, \oo_{X}) \; = \; (\W^2V)^{\otimes \frac{n(n-1)}{2}} \otimes \oo_{X}(1-n).}\leqno{(3.6.3)}
$$
C'est un module inversible sur $X$, et le dual de $\omega^{\otimes 2}$ donne une application
$$
\mu : N^{\otimes 2} \; \rr \; \oo_{X}.
$$

\n {\bf 3.7.\, Proposition }\; {\it Soit $V$ un espace vectoriel de rang 2 sur un corps de caract\'eristique $\neq 2$. Le rev\^etement double $ \pp(V)^n / \mathfrak{A}_{n}\, \rr \,  \pp({\sf TS}^n(V))$ est associ\'e au module inversible 
$$N  = (\W^2V)^{\otimes \frac{n(n-1)}{2}} \otimes \oo_{\pp({\sf TS}^n(V))}(1-n)$$
 et \`a la multiplication $\mu$ pr\'ecis\'ee plus haut.}
\bigskip

\n Reprenons les notations introduites ci-dessus : $Y = \pp(V)^n / \mathfrak{A}_{n}$ et $X =  \pp(V)^n / \mathfrak{S}_{n} = \pp({\sf TS}^n(V))$. L'application $\mu : N^{\otimes 2} \; \rr \; \oo_{X}$ donne 
 une structure de $\oo_{X}$-alg\`ebre sur $\oo_{X} \oplus N$. Consid\'erons la duale $N \rr f_{\star}(\oo_{Y})$ de l'application $\omega$  (3.6.2) ; il faut montrer qu'elle  induit un isomorphisme d'alg\`ebres
$$
\oo_{X} \oplus N \; \rr \; f_{\star}(\oo_{Y}).
$$
Cette d\'emonstration passe par des restrictions \`a des ouverts affines de $X$.
\medskip

\n {\bf 3.7.1}\; Consid\'erons de nouveau le morphisme $\pi :  {\bf P}(V)^n \longrightarrow \; {\bf P}({\sf TS}^n(V)) = X$. Soit $x \in V$ un \'el\'ement non nul ; notons $D(x) \subset \pp(V)$ l'ouvert o\`u $\alpha(x)$ engendre $L$. L'\'el\'ement $x^{\otimes n} \in {\sf TS}^n(V)$ d\'efinit, de m\^eme,  l'ouvert $D(x^{\otimes n}) \subset {\bf P}({\sf TS}^n(V)) = X$ au-dessus duquel l'image $\alpha_{X}(x^{\otimes n})$ de cet \'el\'ement engendre $\oo_{X}(1)$. Par image r\'eciproque, on trouve, compte-tenu de (3.2.2), un isomorphisme
$$
\pi^{\star}(\oo_{X} \stackrel{\alpha_{X}(x^{\otimes n})}{\rr} \oo_{X}(1) ) \; \simeq \; (\oo_{{\bf P}(V)^n} \xrightarrow{\alpha_{1}(x)\otimes \cdots \otimes \alpha_{n}(x)} L_{1}\otimes \cdots \otimes L_{n})
$$
\n Par suite, on a
$$
\pi^{-1}(D(x^{\otimes n}))\;  = \; D(x) \times \cdots \times D(x).
$$
\medskip

Notons aussi que les puissances $x^{\otimes n}$ engendrent ${\sf TS}^n(V)$ d\`es que le corps de base contient $n+1$ \'el\'ements, comme on peut le d\'eduire de (3.2.3), et qu'alors les ouverts de la forme $D(x^{\otimes n})$ recouvrent $X$. La d\'emonstration de {\bf 3.7} autorise les extensions du corps de base. Ainsi, il suffit de faire la d\'emonstration pour les rev\^etements doubles de la forme
$$
D(x)^n / \mathfrak{A}_{n} \; \rr \; D(x^{\otimes n}).
$$
\medskip

\n {\bf 3.7.2}\; Soit donc $\{e_{0}, e_{1}\}$ une base de $V$ ; d\'esignons par $\mathcal{V}$ le sch\'ema affine $D(e_{1})^n  \subset {\bf P}(V)^n$. Pour chaque $i$, $\alpha_{i}(e_{1})$ engendre $L_{i} | \mathcal{V}$, donc il existe $T_{i} \in {\rm H}^0(\mathcal{V}, \oo_{\mathcal{V}})$ tel que 
$$
\alpha_{i}(e_{0}) = T_{i}\alpha(e_{1}).
$$
On a un isomorphisme 
$$
{\rm H}^0(\mathcal{V}, \oo_{\mathcal{V}}) \simeq {\rm H}^0(D(e_{1})^n ) \simeq {\rm H}^0(D(e_{1}))^{\otimes n} \simeq K[T_{1}, \ldots , T_{n}] .
$$
Explicitons de m\^eme l'anneau du sch\'ema affine $\mathcal{U} = D(e_{1}^{\otimes n}) \subset X$. Avec les notations de {\bf 3.2.}, la section inversible s'\'ecrit  $e_{1}^{\otimes n} = {\bf e}_{0}$ ; pour $p = 1, 2, \ldots n$, si on pose $S_{p} = {\bf e}_{p}/{\bf e}_{0} \in {\rm H}^0(\mathcal{U}, \oo_{\mathcal{U}})$, on a un isomorphisme $K[S_{1}, \ldots , S_{n}] \; \simeq \; {\rm H}^0(\mathcal{U}, \oo_{\mathcal{U}})$ ;  le morphisme $\pi : \mathcal{V} \rr \mathcal{U}$ correspond \`a l'inclusion de $K$-alg\`ebres

$$
K[S_{1}, \ldots , S_{n}] \;   \; \subset \;
K[T_{1}, \ldots , T_{n}].
$$
o\`u $S_{p}$ est identifi\'e au polyn\^ome sym\'etrique \'el\'ementaire de degr\'e $p$ en les $T_{i}$.
La fin de cette v\'erification est classique : introduisons le polyn\^ome 
$$
{\sf V}(T_{1}, \ldots, T_{n}) \; = \; \prod_{i<j} (T_{i} - T_{j}).
$$

L'anneau d'invariants $K[T_{1}, \ldots , T_{n}]^{\mathfrak{A}_{n}}$ est un module libre de rang 2 sur $K[S_{1}, \ldots , S_{n}]$, de base $\{1, {\sf V}\}$. Rappelons-en la d\'emonstration : on choisit une transposition $\tau$, de sorte que $\mathfrak{S}_{n} = \mathfrak{A}_{n} \sqcup \tau \mathfrak{A}_{n}$ ; soit $P$  un polyn\^ome invariant sous le groupe altern\'e ;  on \'ecrit
$$
P = \frac{1}{2}(P + {}^{\tau}P) + \frac{1}{2}(P - {}^{\tau}P) = P^+ + P^-.
$$
Comme $P$ est invariant par permutation paire, pour toute transposition $\sigma$, on a \, ${}^{\sigma}P = {}^{\tau \tau \sigma}P = {}^{\tau}({}^{\tau \sigma}P) = {}^{\tau}P$ ; par suite, le premier facteur $P^+$ est sym\'etrique, et on a \, ${}^{\sigma}P^- = - P^-$ ; prenant pour $\sigma$ la transposition $\{i, j\}$, on en d\'eduit que le polyn\^ome $T_{i}-T_{j}$ divise $P^-$, et, de proche en proche, que ${\sf V}$ divise $P^-$ ; si on \'ecrit $P^- = {\sf V}Q$, on constate que le polyn\^ome $Q$ est sym\'etrique. 

Le carr\'e de ${\sf V}$ est le discriminant du polyn\^ome, en $T$, g\'en\'erique ([Bour] A IV 6.7, formule (46))
$$
\prod_{i}(T -T_{i}) \, = \, T^n - S_{1}T^{n-1} + \cdots + (-1)^nS_{n}
$$
C'est un polyn\^ome sym\'etrique, i.e on a  ${\sf V}(T_{1}, \ldots , T_{n})^2  \in K[S_{1}, \ldots , S_{n}]$. D'o\`u la description habituelle de ce  rev\^etement dans le cas polynomial.
\medskip

\n Pour achever la d\'emonstration de ${\bf 3.7}$ il reste \`a relier ${\sf V}$ \`a la section $w$ de $ {\rm H}^0(\mathcal{V}, {\sf M})$ introduite en ${\bf 3.5.}$. Or, l'application (3.5.1) 
$\W^2V \rr L_{i}\otimes L_{j}$ s'\'ecrit ici 
$$ e_{0}\wedge e_{1} \longmapsto T_{i}\alpha_{i}(e_{1})\otimes \alpha_{j}(e_{1}) - T_{j}\alpha_{i}(e_{1})\otimes \alpha_{j}(e_{1}) = (T_{i} - T_{j})\alpha_{i}(e_{1})\otimes \alpha_{j}(e_{1}).$$
Il est alors clair que la section $w : \oo_{\mathcal{V}} \rr {\sf M}| \mathcal{V}$ correspond \`a l'application 
$$
(\W^2V)^{\otimes \frac{n(n-1)}{2}} \rr  (L_{1}\otimes \cdots \otimes L_{n})^{\otimes (n-1)}, \qquad (e_{0}\wedge e_{1})^{\otimes \frac{n(n-1)}{2}} \; \longmapsto {\sf V}.(\alpha_{1}(e_{1})\otimes \cdots \otimes \alpha_{n}(e_{1}))^{\otimes (n-1)}
$$
\bigskip

\n {\bf 3.8.}\; Exemple : le rev\^etement $\pp_{1}\times \pp_{1} \rr \pp_{2}$
\medskip

Appliquant ce qui pr\'ec\`ede lorsque $n = 2$, on trouve que le module inversible associ\'e  \`a ce rev\^etement est isomorphe \`a  $ \oo_{{\bf P}_{2}}(-1)$, et que la multiplication est donn\'ee par la section
$$
T_{1}^2 - 4T_{0}T_{2} \, : \, \oo_{{\bf P}_{2}}(-2)\; \longrightarrow \; \oo_{{\bf P}_{2}} \leqno{(3.8.1)}
$$
(On a identifi\'e $\pp_{2}$ et ${\rm Proj}(K[T_{0}, T_{1}, T_{2}])$). Cela m\'erite d'\^etre reli\'e \`a la formule
$$
(X_{1}Y_{2}\, -\, Y_{1}X_{2})^2\; = \; (X_{1}Y_{2} + Y_{1}X_{2})^2 \, - \, 4(X_{1}X_{2})(Y_{1}Y_{2}).\leqno{(3.8.2)}
$$

\n 
\n En effet,  le morphisme $\pi : \pp_{1}\times \pp_{1} \rr \pp_{2}$ s'\'ecrit en coordonn\'ees homog\`enes :
$$
(x_{1} : y_{1})\, , \, (x_{2} : y_{2}) \; \longmapsto \; (x_{1}x_{2} : x_{1}y_{2}+x_{2}y_{1} : y_{1}y_{2}) 
 $$
 Ainsi, sur $\pp_{1}\times \pp_{1}$, le carr\'e de la section ${X_{1}}\otimes {Y_{2}}\, -\, Y_{1}\otimes X_{2}$ du faisceau inversible  $L_{1}\otimes L_{2} = \pi^{\star}\oo_{{\bf P}_{2}}(1)$, s'exprime par (3.8.2) en fonction des sections invariantes $T_{0} = X_{1}X_{2}, \, T_{1} = X_{1}Y_{2}+Y_{1}X_{2}$  et $T_{2} = Y_{1}Y_{2}$, lesquelles sont des sections de $\oo_{{\bf P}_{2}}(1)$. La propri\'et\'e universelle {\bf 1.5} conduit \`a un morphisme 
 $$
 \pp_{1}\times \pp_{1} \rr {\mathcal Spec}(\oo_{{\bf P}_{2}} \oplus \oo_{{\bf P}_{2}}(-1))
 $$
 D'apr\`es {\bf 3.7}, c'est un isomorphisme.
 
 \n Notons que $\pi$ induit un isomorphisme de la diagonale $\Delta \subset   \pp_{1}\times \pp_{1}$ (dont l'\'equation est  $X_{1}Y_{2}\, -\, Y_{1}X_{2} = 0$)  sur la conique de ${\bf P}_{2}$ d'\'equation  $T_{1}^2 - 4T_{0}T_{2}  = 0$, laquelle est le lieu de diramation de $\pi$.
 \newpage 
%%%%%%%%%%%%%%%%%%%%%%%%%%%%%%%%%%%%%%%%

\noindent {\bf 4.\; Le rev\^etement associ\'e \`a une forme quadratique}
\bigskip

\noindent Soit $E$ un $\OO$-module localement libre de rang deux, et $L$ un $\OO$-module inversible.
On consid\`ere une application quadratique $q : E \rr L$, ou plut\^ot, l'application lin\'eaire qui lui est associ\'ee
$$
\varphi : {\mathcal Sym}^2(E) \quad \longrightarrow \quad L. 
$$
\noindent Il est essentiel de ne pas se limiter aux formes  bilin\'eaires usuelles, qui sont valeurs dans $\OO$.
\medskip

 Ce paragraphe donne la construction, \`a partir de $\varphi$,  d'un rev\^etement de rang 2, $f : Y \rightarrow X$ et d'un module quasi-coh\'erent $\E$ sur $Y$, tels que $f_{\star}(\E) = E$.
 
 \noindent On montre ensuite que si $\varphi$ est surjective, alors $\E$ est  un $\mathcal{O}_{Y}$-module inversible.
\bigskip

\noindent {\bf 4.1.}\; Posons $N = {\mathcal Hom}(L, \mathsf{\Lambda}^2E)$ ; c'est un $\OO$-module inversible (On omet $\OO$ en indice dans $\bigotimes$ et dans ${\mathcal Hom}$ car ici tout prend place dans la cat\'egorie de $\OO$-modules). Partant de la forme $\varphi$, comme ci-dessus, on d\'efinit une application $\OO$-lin\'eaire
$$
u : N\otimes E \; \longrightarrow \; E \leqno{(4.1.1)}
$$
en composant les applications suivantes
$$
N\otimes E \longrightarrow N\otimes {\mathcal Hom}(E, L)\quad  \widetilde{\longrightarrow} \quad {\mathcal Hom}(E, \wedge^2E)  \quad \widetilde{\longleftarrow}\quad E,
$$
o\`u celle de gauche est associ\'ee \`a $\varphi$, et avec le facteur $\frac{1}{2}$. Plus pr\'ecisement, et en termes de sections locales $x, y$ de $E$, et d'une  section locale $\alpha$ de $N= {\mathcal Hom}(L, \W^2E)$, on d\'efinit $u$ par l'\'egalit\'e
$$
\boxed{x \wedge u(\alpha \otimes y)\; = \; \frac{1}{2}\alpha(\varphi(x.y))}\leqno{(4.1.2)}
$$

\n Il y a deux isomorphismes\,  \og naturels\fg \, entre $E$ et ${\mathcal Hom}(E, \wedge^2E) $ ; ils diff\`erent par un signe ; j'ai  choisi  l'isomorphisme
$E \rightarrow {\mathcal Hom}(E, \wedge^2E)$ donn\'e par $y \longmapsto (x\mapsto x\wedge y)$ ; il semble faciliter l'\'ecriture de plusieurs formules. Le coefficient $\frac{1}{2}$ est n\'ecessaire pour que $x \mapsto x \wedge u(\alpha \otimes x)$ soit li\'e \`a la forme quadratique associ\'ee \`a $\varphi$, $x\mapsto \frac{1}{2}\varphi(x^2)$. On trouvera en {\bf 4.2.} une traduction matricielle de $u$.

\medskip

\noindent On \'ecrira  souvent $\alpha x$ pour $u(\alpha \otimes x)$ ; l'application $x \mapsto \alpha x$ est un endomorphisme de $E$ not\'e $\alpha_{E}$. 
\medskip

La construction de $u$ \`a partir de $\varphi$ commute aux changements de base puisque $E, L$ et $N$ sont localement libres ; en consid\'erant, pour $x \in X$,  le changement de base ${\rm Spec}(\kappa(x)) \rightarrow X$, on voit que si $\varphi$ est non nulle en $x$, alors l'application $E\otimes \kappa(x) \rr {\mathcal Hom}(E, L)\otimes \kappa(x)$ est non nulle, et par suite $1_{\kappa(x)}\otimes u \neq 0$.
\bigskip

On d\'efinit une application lin\'eaire
$$
\mu
: N^{\otimes 2}\; \longrightarrow \; \OO,  
$$
en posant 
$$
\mu(\alpha^{\otimes 2}) = - \det(\alpha_{E})\leqno{(4.1.3)}
$$

\n On peut aussi caract\'eriser $\mu$ par la commutativit\'e du carr\'e suivant 
$$
\begin{CD}
\W^2(N\otimes E) @>{\wedge^2u}>> \W^2E\\
@V\wr VV @VV= V\\
N^{\otimes 2}\otimes \W^2E @>>{-\mu \otimes 1}> \W^2E
\end{CD}
$$
c'est-\`a-dire par la formule
$$
\mu(\alpha \otimes \beta) x\wedge y \, = \,- \, \alpha_{E}(x) \wedge \beta_{E}(y).\leqno{(4.1.4)}
$$

\noindent Cette application $\mu$ permet de munir le $\OO$-module de rang 2, 
$$
\A \; = \; \OO \oplus N
$$
d'une structure de $\OO$-alg\`ebre.  Le morphisme 
$$f : Y= {\rm Spec}(\A) \longrightarrow X
$$
 est donc fini, localement libre de rang deux. 
 
 Le ferm\'e de diramation ({\bf 1.3})  de $f$ est d\'efini par l'id\'eal $
 \mathcal {I} = {\rm Im}(N^{\otimes 2} \stackrel{\mu}{\rightarrow} \OO) $ ; par suite, $Y \rr X$ est \'etale si et seulement si $u$ est un isomorphisme ; on dit alors parfois que $q$ est non d\'eg\'en\'er\'ee, ou non singuli\`ere.
\bigskip

\n {\bf 4.2.}\; Expression locale
\medskip

Lorsque $X$ est affine d'anneau $R$, et que $E$ et $L$ sont des modules libres, ces constructions ont la traduction matricielle suivante.

Soit $\{e_{1}, e_{2}\}$ une base de $E$, de sorte que $\{e_{1}^2, e_{1}e_{2}, e_{2}^2\}$ est une base de ${\rm Sym}^2(E)$ ; le choix d'une base $\{\varepsilon \}$ de $L$ permet d'\'ecrire $\varphi = \psi \epsilon$, o\`u  $ \psi $ est une forme lin\'eaire sur ${\rm Sym}^2(E)$. Posons
$$
a = \psi(e_{1}^2), \qquad b = \psi(e_{1}e_{2}), \qquad c = \psi(e_{2}^2) .
$$
On retrouve l'expression usuelle :
$$
\psi((xe_{1}+ye_{2})^2) \; = \; ax^2 + 2bxy + cy^2.
$$
Soit $\alpha$ la base de $N = {\rm Hom}(L, \W^2E)$ d\'efinie par $\alpha(\varepsilon) = e_{1} \wedge e_{2}$. Explicitons la matrice de l'endomorphisme  $\alpha_{E} : E \rr E,\; x \mapsto \alpha x = u(\alpha \otimes x)$. La relation (4.1.2) donne
$$
e_{i}\wedge \alpha e_{j} \; = \; \frac{1}{2}\alpha(\varphi(e_{i}e_{j})) = \frac{1}{2}\psi(e_{i}e_{j})e_{1}\wedge e_{2}.
$$
On en tire  la matrice de $\alpha_{E}$ relativement \`a la base $\{e_{1}, e_{2}\}$ :
$$
 \frac{1}{2} \begin{pmatrix}
\; -b&-c\\
\; a&\; b
\end{pmatrix} \leqno{(4.2.1)}
$$
De m\^eme, on trouve pour l'application $\mu$ de (4.1.3)
$$
\mu(\alpha^{\otimes 2}) \; = \; \frac{1}{4} (b^2 - ac).
$$
Finalement, on voit que l'alg\`ebre $\A$ associ\'ee \`a $\varphi$ est isomorphe \`a
$$
R[T]/(T^2 - ( b^2 -ac)),
$$
o\`u la classe de $T$ correspond \`a la base $2\alpha$ de $N$.
\bigskip

\noindent {\bf 4.3.}\; Structure de $\A$-module sur $E$.
\medskip

\noindent Revenons \`a la situation g\'en\'erale de {\bf 4.1.} On munit $E$ d'une structure de $\A$-module en utilisant l'application $u : N\otimes E \rightarrow E$, laquelle s'\'etend en une application $\OO$-lin\'eaire $\A \otimes E \rightarrow E$. L'associativit\'e (i.e l'\'egalit\'e $\alpha(\beta x) = (\alpha\beta) x$) est cons\'equence de la relation $\alpha(\alpha x) = (\alpha^2)x$, puisque $N$ est un module inversible.

Cette relation  d\'ecoule des \'egalit\'es suivantes :\\

$y \wedge \alpha(\alpha x)\; = \; \frac{1}{2}\alpha(\varphi(y.\alpha x))$,  \hspace{12mm}  d'apr\`es (4.1.2)

\hspace{2,1cm} $ =  \frac{1}{2}\alpha(\varphi(\alpha x . y))$, \hspace{8mm} puisqu'il s'agit du produit dans ${\rm Sym^2}$

\hspace{2,1cm} $ =  \alpha x \wedge \alpha y$, \hspace{15mm} d'apr\`es (4.1.2)

\hspace{2,1cm} $= \det(\alpha_{E}) x\wedge y, $ \hspace{8mm} par d\'efinition de $\det(\alpha_{E})$

\hspace{2,1cm} $ = \alpha^2\, y \wedge x,$ \hspace{15mm}  d'apr\`es (4.1.3)

\hspace{2,1cm} $ = y \wedge (\alpha^2)x,$ \hspace{13mm}  puisque $\alpha^2 \in \OO$.\\

\n Comme cette \'egalit\'e est vraie pour toute section $y$, on voit que l'on a l'\'egalit\'e annonc\'ee $\alpha(\alpha x) = (\alpha^2)x$.
\bigskip

\n {\bf Remarque}\; L'\'egalit\'e (4.1.2) est la cl\'e des constructions du texte : comme $E$ est de rang deux, la connaissance de $x \wedge u(\alpha \otimes y)$ pour tout $x$, d\'etermine $u(\alpha \otimes y)$ ; ainsi, $u$ et $\varphi$ se cod\'eterminent l'une l'autre ; autrement dit, se donner la structure de $\A$-module sur $E$ (c'est-\`a-dire $u$) revient \`a se donner la forme quadratique sur $E$ (c'est-\`a-dire $\varphi$) ; cette formule montre  enfin le r\^ole de l'inversibilit\'e de 2. Voir {\bf 5.4.} pour un \'enonc\'e pr\'ecis.
\bigskip

\noindent {\bf 4.4.}\; Montrons que la surjectivit\'e de $\varphi$ entra\^ine que $E$ est un $\A$-module inversible, {\it i.e.}
que $\E$ est un $\oo_{Y}$-module inversible. 
\medskip

\noindent Pour ce faire, on peut supposer que $X$ est le spectre d'un anneau local $R$, et donc que $E$, $L$ et $N$ sont des $R$-modules libres. Choisissons des bases et reprenons les notations de {\bf 4.2.}\, Comme  le $R$-module ${\rm Sym}^2(E)$ est engendr\'e par les carr\'es ({\bf 0.2}), la surjectivit\'e de $\varphi $ se traduit  par l'existence d'un  $x \in E$ tel que l'\'el\'ement  $\varphi(x.x) = \psi(x.x) \varepsilon $ soit un g\'en\'erateur de $L$, {\it i.e} que $\psi(x.x)$ soit inversible dans $R$. Notant $\alpha$ la base de $N$ telle que $\alpha(\varepsilon) = e_{1} \wedge e_{2}$, la relation
$$
x \wedge \alpha  x \; = \; \frac{1}{2}\psi (x.x) e_{1}\wedge e_{2}
$$
montre que $\{x, \, \alpha  x\}$ est une base du $R$-module $E$, donc que $x$ est un g\'en\'erateur de ce module sur l'anneau $\A = R + N$ ; c'en est m\^eme une base puisque la surjection obtenue $\A \rightarrow E$ concerne deux $R$-modules libres de m\^eme rang, et est donc un isomorphisme.
\vspace{1cm}

%%%%%%%%%%%%%%%%%%
\noindent {\bf 5.\; Forme quadratique sur l'image directe d'un inversible}
\bigskip

Soit $f : Y \rightarrow X$ un rev\^etement de rang deux, et $\E$ un faisceau inversible sur $Y$. Le  $\OO$-module 
$$
E \, =\, f_{\star}(\E)
$$
est donc localement libre de rang deux.

On rappelle d'abord comment munir  canoniquement $E$ d'une forme quadratique \`a valeurs dans le module inversible  $L = {\sf N}_{Y/X}(\E)$, et on montre  que le rev\^etement de rang 2 associ\'e \`a cette forme selon le \S 4 est isomorphe \`a \, $Y \rightarrow X$. La construction du \S 4 est donc r\'eversible,  et cela est pr\'ecis\'e en une \'equivalence de cat\'egories.
\bigskip

\noindent {\bf 5.1.}\; Gardons les notations introduites au \S 1 : 
$
f_{\star}(\mathcal{O}_{Y}) \; =\; \A \; = \; \OO \oplus N, 
$
o\`u $N = {\rm Ker}(\A \stackrel{{\rm Tr}}{\longrightarrow} \OO)$.

\n Le foncteur ${\sf N}_{Y/X}$ associe \`a tout $\oo_{Y}$-module inversible $\E$ le $\OO$-module inversible 
$$
{\sf N}_{Y/X}(\E) \, =\, {\mathcal Hom}_{\OO}(\W^2f_{\star}(\mathcal{O}_{Y}), \, \W^2 f_{\star}(\E)) = {\mathcal Hom}_{\OO}(\W^2 \A, \, \W^2 E)
$$
L'application quadratique $\nu_{\E}$ est \emph{l'application normique universelle}  $ E = f_{\star}(\E) \rr {\sf N}_{Y/X}(\E)$ (voir, pour plus de d\'etails, ([Fer 1], \S 3.3)) ; elle sera  nomm\'ee ici  plus simplement  la \og norme \fg\, (avec des guillemets) ; c'est l'application fournie par le carr\'e ext\'erieur
$$
E  = {\mathcal Hom}_{\A}(\A, E) \; \subset \; {\mathcal Hom}_{\OO}(\A, E)\quad \stackrel{\wedge^2}{\longrightarrow}\quad {\mathcal Hom}_{\OO}(\W^2 \A, \W^2E)
$$
Tout \'el\'ement de $\W^2\A = \W^2(\OO \oplus N)$ s'\'ecrit localement de fa\c{c}on unique sous la forme $1\wedge \alpha$, avec $\alpha$ une secction de $N$, et on a 
$$
\nu_{\E}(x) \, = \, (1\wedge \alpha \mapsto x\, \wedge \, \alpha x). \leqno{(5.1.1)}
$$

\noindent C'est une application polynomiale de degr\'e 2, qui s'\'etend  en l'application lin\'eaire
$$
\varphi : {\mathcal Sym}^2(E) \quad \longrightarrow \quad L,
$$
d\'efinie par $ \varphi(xy) = \nu(x+y)-\nu(x) - \nu(y)$. On trouve
$$
\varphi(xy) = (1 \wedge \alpha \; \longmapsto \; (x\wedge \alpha y  + y \wedge \alpha x))
$$
En fait, les deux termes  dans la seconde parenth\`ese sont \'egaux.
\bigskip

\noindent {\bf 5.1.2.\; Lemme }\; {\it Soit $\alpha$ un endomorphisme d'un module $E$ localement libre de rang deux. Si ${\rm Tr}_{E}(\alpha) = 0$, alors, pour tous $x, y \in E$, on a dans $\W^2 E$,}
$$
x \wedge \alpha y \; = \; y \wedge \alpha x.
$$
\medskip

En effet, introduisons une ind\'etermin\'ee $T$, et calculons de deux fa\c{c}ons $(T-\alpha)(x) \wedge (T-\alpha)(y)$ : si on d\'eveloppe on trouve 
$$
T^2 x\wedge y - T(\alpha x\wedge y + x\wedge \alpha y) + \det(\alpha)x\wedge y,
$$
Par ailleurs, par d\'efinition du d\'eterminant, cet \'el\'ement est aussi \'egal \`a
$$
(T^2 -{\rm Tr}(\alpha)T + \det(\alpha)) x \wedge y
$$
 L'hypoth\`ese ${\rm Tr}(\alpha) = 0$ entra\^ine donc l'\'egalit\'e annonc\'ee. $\Box$
 \bigskip
 
 \n Si $\alpha$ est dans $N = {\rm Ker(Tr}_{\A})$, on a aussi ${\rm Tr}_{E}(\alpha) = 0$ puisque $E$ est un $\A$-module inversible. Finalement, on obtient l'expression suivante pour l'application 
 $$
 \varphi : {\mathcal Sym}^2(E) \quad \longrightarrow \quad L={\mathcal Hom}_{\OO}(\W^2\A, \, \W^2 E)
 $$ 
 $$
 \varphi(xy) \; =\;(1\wedge \alpha \longmapsto  2 x\wedge \alpha y).\leqno{(5.1.3)}
 $$
 \medskip
 
La surjectivit\'e de cette application est une propri\'et\'e locale ; pour la v\'erifier, on peut donc supposer que $\E = \mathcal{O}_{Y}$, et, par suite, que $E = \A = \OO \oplus N$,  et $L =\OO$ ; l'image de $\varphi$ est un id\'eal, et il contient l'\'el\'ement inversible $2 = \varphi(1.1)$.
\bigskip 

\noindent {\bf 5.2.}\; Montrons que le rev\^etement associ\'e \`a $\varphi$ est isomorphe \`a $f$.
\medskip

 Reprenons la d\'emarche du \S 4 : la construction du rev\^etement associ\'e \`a $\varphi$ repose sur le module inversible ${\mathcal Hom}(L, \W^2E)$, not\'e $N'$, et sur l'application (4.1.1)\, $u :N'\otimes E \rr E$. 

\n On d\'efinit un isomorphisme $N \simeq N'$ en associant \`a $\alpha \in N$ l'\'el\'ement $\alpha' \in N' =  {\mathcal Hom}(L, \W^2E)$ d\'efini par
$$
\alpha'(\lambda) = \lambda(1 \wedge \alpha)  \leqno{(5.2.1)}
$$
(Cela a bien un sens puisque $\lambda \in L={\mathcal Hom}_{\OO}(\W^2\A, \, \W^2 E)$ (5.1.3)).

\n Il faut d'abord  v\'erifier que pour tout $\alpha \in N$ et $y \in E$, on a
$$
u(\alpha' \otimes y) = \alpha y \leqno{(5.2.2)}
$$
o\`u le produit  $\alpha y$  est donn\'e par la structure initiale de $\A$-module sur $E$. Or, pour tout $x \in E$, on a
$$
x \wedge u(\alpha'\otimes y)\; \stackrel{(4.1.2)}{=}\;  \frac{1}{2}\alpha'(\varphi(xy)) \; \stackrel{(5.2.1)}{=}\;\frac{1}{2}\varphi(xy)(1\wedge \alpha)\;  \stackrel{(5.1.3)}{=}\;  x \wedge \alpha y.
$$
D'o\`u l'\'egalit\'e (5.2.2). Il faut ensuite v\'erifier que l'isomorphisme lin\'eaire $\A = \OO \oplus N \, \simeq \, \OO \oplus N'$ est un isomorphisme d'alg\`ebres, c'est-\`a-dire que le carr\'e dans $\A$ d'un \'el\'ement $\alpha \in N$, soit $ - {\rm norm}_{\A}(\alpha)$, est \'egal au carr\'e de $\alpha'$  soit $- \det(\alpha'_{E})$ (4.1.3). Or, l'\'egalit\'e (5.2.2) dit que l'endomorphisme $\alpha'_{E} = (y \mapsto u(\alpha'\otimes y))$ est \'egal \`a $(y \mapsto \alpha y)$ ; comme $E$ est un $\A$-module inversible, le d\'eterminant de ce dernier est \'egal \`a ${\rm norm}_{\A}(\alpha)$ ; d'o\`u la compatibilit\'e  aux produits, qui \'etait  annonc\'ee.
\bigskip

\n {\bf 5.3}\; Fonctorialit\'e des ces constructions.
\medskip

Au couple $(Y, \E)$ form\'e d'un rev\^etement double de $X$ et d'un $\oo_{Y}$-module inversible $\E$, on associe donc une forme quadratique partout non nulle sur $X$,\, $q : E \rr L$. Il s'agit ici de d\'egager des morphismes de couples qui  induisent des morphismes de formes quadratiques. Notons d'abord qu'une restriction \'evidente s'impose : on ne peut consid\'erer que les morphismes de rev\^etements  de $X$,
$$
\xymatrix{
Y' \ar[dr]_{f' }\ar[rr]^g && Y \ar[dl]^f \\
&X&
 }
$$
qui sont, au moins,  compatibles aux normes usuelles, au sens o\`u le diagramme suivant doit \^etre  commutatif
$$
\xymatrix{
f_{\star}\oo_{Y} \ar[dr]_{{\rm norm}}\ar[rr]^\psi &&f'_{\star}\oo_{Y'} \ar[dl]^{\rm norm} \\
&\OO&
 }
$$
Le morphisme $\psi$ d\'esigne l'image directe par $f$ du morphisme canonique $\oo_{Y} \rr g_{\star}\oo_{Y'}$, et cette commutativit\'e est requise pour tous les triangles obtenus par changement de base sur $X$ ; autrement dit, les normes doivent \^etre consid\'er\'ees comme des \emph{lois polyn\^omes} et non comme de simples applications. Introduisant une ind\'etermin\'ee $T$, on doit donc avoir, pour toute section locale $y$ de $f_{\star}\oo_{Y} $, 
$$
{\rm norm}_{Y/X}(T- y)  \; =\; {\rm norm}_{Y'/X}(T- \psi(y)) 
$$
c'est-\`a-dire
$$
T^2 - {\rm Tr}_{Y/X}(y) T + {\rm norm}_{Y/X}(y) \; =\; T^2 - {\rm Tr}_{Y'/X}(\psi(y)) T + {\rm norm}_{Y'/X}(\psi(y))
$$
Cette \'egalit\'e entre polyn\^omes implique la suivante 
$$
{\rm Tr}_{Y/X}(y) \; = \; {\rm Tr}_{Y'/X}(\psi(y)) .
$$
On dira donc parfois que $g$ est \emph{compatible aux traces}.
En utilisant les d\'ecompositions $f_{\star}\oo_{Y} = \OO \oplus N$  et  $f'_{\star}\oo_{Y'} = \OO \oplus N'$, on voit que la compatibilit\'e aux traces s'\'ecrit finalement
$$
\psi(N) \; \subset \; N' . \leqno{(5.3.1)}
$$
En fait on n'obtiendra de r\'esultats satisfaisants que sous une hyppoth\`ese un peu plus forte (cf. ({\bf 1.6})), \`a savoir :
$$
\psi \quad {\rm est\; injective} .\leqno{(5.3.2)}
$$
\medskip

\n {\bf 5.3.3}\; {\it Gardons les notations introduites, et consid\'erons un morphisme $g : Y' \rr Y$ de rev\^etements doubles de $X$ ; soit $\E$ un module inversible sur $Y$, et $\E' = g^{\star}\E$ son image r\'eciproque sur $Y'$. On suppose que l'homomorphisme $\psi : f_{\star}\oo_{Y} \rr  f'_{\star}\oo_{Y'}$ associ\'e \`a $g$ est \emph{injectif}. Alors $g$ induit un morphisme des \og normes \fg, au sens suivant : il existe un isomorphisme $\omega : {\sf N}_{Y/X}(\E) \simeq {\sf N}_{Y'/X}(\E')$ rendant commutatif le diagramme }
$$
\begin{CD}
f_{\star}\E @>{\theta}>> f'_{\star}\E'\\
@V{\nu_{\E}}VV @VV{\nu_{\E'}}V\\
{\sf N}_{Y/X}(\E) @>>{\omega}> {\sf N}_{Y'/X}(\E') ,
\end{CD}
$$
{\it o\`u $\theta$ d\'esigne l'image directe par $f$ de l'application canonique $\E \rr g_{\star}g^{\star}\E$.}
\bigskip

All\'egeons, comme plus haut, les notations en posant $\A = f_{\star}\oo_{Y},\, \A' = f'_{\star}\oo_{Y'}, \, E = f_{\star}\E$ et $E' = f'_{\star}\E'$. Pour d\'efinir $\omega$, introduisons les deux applications \'evidentes
$$
{\sf N}_{Y/X}(\E) = {\mathcal Hom}(\W^2\A, \W^2E) \; \xrightarrow{(1, \wedge^2\theta)} \;{\mathcal Hom}(\W^2\A, \W^2E')\; \xleftarrow{(\wedge^2\psi, 1)} \;{\mathcal Hom}(\W^2\A', \W^2E') = {\sf N}_{Y'/X}(\E')
$$
Il s'agit de v\'erifier qu'elles sont injectives et qu'elles ont m\^eme image. Ce sont l\`a des propri\'et\'es locales sur $X$, si bien qu'on peut supposer, $\E$ \'etant localement isomorphe \`a $\oo_{Y}$, que $E = \A$, $E' = \A'$ et $\theta = \psi$. L'injectivit\'e de $\wedge^2\psi$ provient de l'injectivit\'e de $\psi$, et celle de sa duale provient du fait que ${\rm Coker}(\wedge^2\psi)$ est un module de torsion. Par ailleurs, les deux applications
$$
{\mathcal Hom}(\W^2\A, \W^2\A) \; \xrightarrow{(1, \wedge^2\psi)} \;{\mathcal Hom}(\W^2\A, \W^2\A')\; \xleftarrow{(\wedge^2\psi, 1)} \;{\mathcal Hom}(\W^2\A', \W^2\A')
$$
ont la m\^eme image : le sous-module engendr\'e par $\wedge^2\psi \in {\rm Hom}(\W^2\A, \W^2\A')$.

\n La compatibilit\'e aux \og normes\fg \, se voit donc par la commutativit\'e du diagramme suivant
$$
\xymatrix{
E\ar[d]_{\nu} \ar[rr]^{\theta} && E' \ar[d]^{\nu}\\
{\mathcal Hom}(\W^2\A, \W^2E)\ar[dr]_{(1, \wedge^2\theta)}&& {\mathcal Hom}(\W^2\A', \W^2E')\ar[dl]^{(\wedge^2\psi, 1)} \\
&{\mathcal Hom}(\W^2\A, \W^2E')&
 } 
$$
laquelle se v\'erifie en suivant le destin d'une section $x$ de $E$ :
$$
\xymatrix{
x\ar@{|->}[d] \ar@{|->}[rr] && 1\otimes x \ar@{|->}[d]\\
(1\wedge \alpha \mapsto x\wedge \alpha x)\ar@{|->}[dr]&& (1\wedge \alpha' \mapsto 1\otimes x\, \wedge \,  \alpha'\otimes x) \ar@{|->}[dl]\\
&(1\wedge \alpha \mapsto 1\otimes x \wedge 1\otimes \alpha x)&
 }
$$
(Remarquer que, dans $E' = \A'\otimes_{\A}E$, on a $\psi(\alpha)\otimes x = 1\otimes \alpha x$ pour $\alpha \in N$).
\medskip

Examinons maintenant la fonctorialit\'e du passage d'une forme quadratique \`a un rev\^etement double, expliqu\'e au \S 4.
\bigskip

\n {\bf 5.3.4. }\; {\it Soient $E$ et $E'$ des $\OO$-modules localement libres de rang 2. Soient $q : E \rr L$ et $q' : E' \rr L$ deux applications quadratiques, partout non nulles, \`a valeurs dans un m\^eme $\OO$-module inversible $L$. Soit $f :Y \rr X$ et $f' : Y' \rr X$ les rev\^etements doubles associ\'es \`a $q$ et $q'$.}

\n {\it Alors, \`a toute application lin\'eaire \emph{injective} $\theta : E \rr E'$ telle que $q = q' \theta$, correspond un morphisme de rev\^etements $g : Y' \rr Y$ dont le morphisme d'alg\`ebres associ\'e}
 $$
 \psi : f_{\star}\oo_{Y}= \OO \oplus N \; \rr \; f'_{\star}\oo_{Y'} = \OO \oplus N'
 $$
 {\it est injectif (il  v\'erifie donc la relation $\psi(N) \subset N'$).}
\bigskip

La d\'efinition de la restriction de $\psi$ \`a $N$ s'impose d'elle-m\^eme : c'est l'application
$$
N = {\mathcal Hom}(L, \W^2E) \xrightarrow{(1, \W^2\theta)} {\mathcal Hom}(L, \W^2E') = N'
$$
Il faut v\'erifier tout ce qui est sous-entendu dans le mot \emph{correspond} de la conclusion. \`A savoir les points {\it a)}, {\it b)}  et {\it c)} suivants .
\medskip

\n {\it a)}\; L'op\'eration de $N$ sur $E$ est transport\'ee par $\psi$ en l'op\'eration de $N'$ sur $E'$ ; autrement dit, en utilisant les applications $u$ de (4.1.1), le carr\'e suivant est commutatif
$$
\begin{CD}
N \otimes E @>{\psi \otimes \theta}>> N'\otimes E'\\
@VuVV  @VVu'V\\
E @>>\theta > E'
\end{CD}
$$
Au vu de la relation (4.1.2), il faut v\'erifier que pour toutes sections locales  $y \in E$, $x' \in E'$ et  $\alpha \in N$, on a
$$
x' \wedge u'(\psi(\alpha)\otimes \theta(y)) = x' \wedge \theta(u(\alpha\otimes y)). \leqno{(5.3.4.1)}
$$
Dans le cas o\`u $x'$ provient de $E$, c'est-\`a-dire si $x' = \theta(x)$, cette \'egalit\'e d\'ecoule de la d\'efinition de $u$ (on note $\varphi$ et $\varphi'$ les formes bilin\'eaires associ\'ees \`a $q$ et $q'$ respectivement); en effet, on a les \'egalit\'es suivantes :
\medskip

\n $\theta(x) \wedge u'(\psi(\alpha)\otimes \theta(y)) = \frac{1}{2}\psi(\alpha)(\varphi'(\theta(x)\theta(y)))$, \hspace{4cm} par d\'efinition (4.1.2)

\n \hspace{3,4cm} $= \frac{1}{2}\psi(\alpha)(\varphi(xy))$, \hspace{5cm} puisque $q' \theta = q$

\n \hspace{3,4cm} $= \frac{1}{2} \wedge^2\theta (\alpha(\varphi(xy)))$, \hspace{4,5cm} par d\'efinition de $\psi$

\n \hspace{3,4cm} $= \theta(x) \wedge \theta(u(\alpha\otimes y))$, \hspace{4,2cm} (4.1.2).
\medskip

Pour un $x'$ g\'en\'eral, on se ram\`ene au cas pr\'ec\'edent par la remarque suivante. L'\'egalit\'e (5.3.4.1) peut se v\'erifier localement sur $X$, ce qui permet de supposer que $E = E'$, et que $\theta$  est donc un endomorphisme ; son d\'eterminant est une section r\'eguli\`ere puisque $\theta$ est suppos\'ee injective. Soit $t$ cette section. Comme $\W^2E'$ est inversible, $t$ est r\'eguli\`ere pour ce module, et il suffit donc de v\'erifier l'\'egalit\'e $tx' \wedge u'(\psi(\alpha)\otimes \theta(y)) = tx' \wedge \theta(u(\alpha\otimes y))$. L'existence du \emph{cotranspos\'e} $\tilde{\theta}$ de $\theta$, et la relation $\theta\circ \tilde{\theta} = t_{E'}$ montrent  que $tE' \subset \theta(E)$, ce qui r\'eduit la v\'erification au cas d\'ej\`a trait\'e.
\medskip

\n {\it b)}\; L'application $\psi$ est compatible aux multiplications {\it i.e.} pour $\alpha \in N$, on a
$$
\mu'(\psi(\alpha)^{\otimes 2}) = \mu(\alpha^{\otimes 2})
$$
Cela d\'ecoule de la formule (4.1.3)  $\mu(\alpha^{\otimes 2}) = - \det(\alpha_{E})$, et de {\it a)}.
Cela entra\^ine que $\psi$ s'\'etend en un morphisme de $\OO$-alg\`ebres $\A \rr \A'$, d'o\`u le morphisme de rev\^etements  $g : Y' \rr Y$.
\medskip

\n {\it c)}\; Soient $\E$ et $\E'$ les modules inversibles associ\'es respectivement \`a $E$ et $E'$. Le point {\it a)} implique que $\theta$ induit une application de $\oo_{Y'}$-modules inversibles $\chi : g^{\star}\E \rr \E'$. C'est un isomorphisme.

\n En effet, d'apr\`es le r\'esultat direct ({\bf 5.3.3}) on dispose d'un isomorphisme $\omega :  {\sf N}_{Y/X}(\E) \simeq {\sf N}_{Y'/X}(g^{\star}\E)$ compatible aux normes. Mais, d'apr\`es ({\bf 5.2}), les deux modules inversibles  ${\sf N}_{Y/X}(\E)$ et ${\sf N}_{Y'/X}(\E')$ sont canoniquement isomorphes \`a $L$ ; on en d\'eduit que ${\sf N}_{Y'/X}(\chi)$ est un isomorphisme (se souvenir que $q$ et $q'$ sont partout non nulles). En explicitant  ${\sf N}_{Y'/X}(g^{\star}\E)$ et ${\sf N}_{Y'/X}(\E')$, on voit l'isomorphisme
$$
{\sf N}_{Y'/X}(\chi) : {\mathcal Hom}(\W^2\A' , \W^2 \A'\otimes_{\A}E) \rr    {\mathcal Hom}(\W^2\A', \W^2E')
$$
Il montre que  l'application $\W^2 \A'\otimes_{\A}E \rr    \W^2E'$ est un isomorphisme, donc que $\chi$ lui-m\^eme en est un.

Cela ach\`eve la d\'emonstration de {\bf 5.3.4.}
\vspace{1cm}

\n {\bf 5.4.\; R\'esum\'e}
\medskip

Soit $X$ un sch\'ema o\`u $2$ est inversible.
\medskip

D\'esignons par $\mathcal{R}$ la cat\'egorie dont les objets sont les couples $(Y, \E)$, o\`u $f :Y\rightarrow X$ est un rev\^etement double de $X$, et $\E$ est un $\oo_{Y}$-module inversible. Une fl\`eche dans $\mathcal{R}$ de $(Y, \E)$ vers $(Y', \E')$ est constitu\'ee d'un morphisme de rev\^etements $g : Y' \rr Y$ et d'un isomorphisme $g^{\star}\E \simeq \E'$ ; on impose de plus \`a $g$ la condition que le morphisme induit $\psi : f_{\star}\oo_{Y} \rr f'_{\star}\oo_{Y'}$ soit \emph{injectif}.
\medskip

D\'esignons par $\mathcal{Q}$ la cat\'egorie dont les objets sont les applications quadratiques partout non nulles  $q : E \rr L$, de source un $\OO$-module localement libre de rang 2, \`a valeurs dans un inversible. Une fl\`eche dans $\mathcal{Q}$ de $q : E \rr L$ vers  $q' : E' \rr L'$ est constitu\'ee d'une application lin\'eaire \emph{injective}  $\theta : E \rr E'$ et d'un isomorphisme $\omega : L \rr L'$ tels que le diagramme suivant soit commutatif
$$
\begin{CD}
E @>\theta >> E'\\
@VqVV @VVq'V\\
L @>>\omega> L'
\end{CD}
$$
\medskip

Alors les cat\'egories $\mathcal{R}$  et  $\mathcal{Q}$ sont \'equivalentes.
\medskip

Plus pr\'ecisement, l'application sur les objets
$$
{\sf N} : \mathcal{R} \; \rr \;\mathcal{Q},\qquad (Y, \E) \longmapsto (\nu_{\E} : f_{\star}\E \rr {\sf N}_{Y/X}(\E))
$$
se prolonge en un foncteur covariant ({\bf 5.3.3}). Et il existe  un foncteur covariant
$$
{\sf A} : \mathcal{Q} \; \rr \;\mathcal{R}
$$
d\'efini sur les objets en {\bf 4.1}, et sur les fl\`eches en {\bf 5.3.4}, tel que les foncteurs ${\sf N}\circ{\sf A}$  et  ${\sf A}\circ{\sf N}$ soient isomorphes aux foncteurs identit\'e.

\vspace{2cm}

%%%%%%%%%%%%%%%%%%%%%%%%%%%%%%%%
\noindent {\bf 6. \; Formes bilin\'eaires sym\'etriques  et polyn\^omes homog\`enes de degr\'e deux}
\\

 \noindent Ces deux notions sont duales l'une de l'autre et on peut, le plus souvent, les identifier sans obscurcir le propos ; mais ici, il faut les distinguer et pr\'eciser comment on passe de l'une \`a l'autre.\\
 
 \noindent Soit $E$ un $\OO$-module localement libre de rang deux, et $L$ et $M$ des $\OO$-modules inversibles.
 
\noindent  
Un polyn\^ome homog\`ene de degr\'e 2 sur $E$ (on devrait pr\'eciser : \og tordu par $M$\fg\,) est une application lin\'eaire 
$$
\gamma: \OO \quad  \longrightarrow \quad {\mathcal Sym}^2(E) \otimes_{\OO}M.
$$
On va rappeler comment il lui correspond une application lin\'eaire 
$$
\varphi : {\mathcal Sym}^2(E) \quad \longrightarrow \quad L.
$$

\noindent Pour passer d'une notion \`a l'autre, on utilise le
\\

{\bf 6.1.\, Lemme }\; {\it Soit $X$ un sch\'ema sur lequel 2 est inversible, et $E$ un $\OO$-module localement libre de rang deux. Alors l'application lin\'eaire}
$$
E^{\otimes 4} \; \longrightarrow \; (\W^2E)^{\otimes 2}, \qquad x_{1}\otimes x_{2}\otimes x_{3}\otimes x_{4} \; \longmapsto \; \, (x_{1}\wedge x_{3}) \otimes (x_{2}\wedge x_{4}) + (x_{1}\wedge x_{4}) \otimes (x_{2}\wedge x_{3}\,)
$$
{\it est invariante si on permute 1 et 2, et si on permute 3 et 4, de sorte qu'elle passe aux quotients et d\'efinit une application 
}$$
\sy^2(E) \otimes \sy^2(E) \; \longrightarrow \; (\W^2E)^{\otimes 2},
$$
{\it d'o\`u, finalement, une application}
$$
\sy^2(E) \; \longrightarrow \; {\mathcal Hom}_{\OO}(\sy^2(E), \,(\W^2E)^{\otimes 2}).
$$
{\it C'est un isomorphisme.}
\vspace{7mm}

On le v\'erifie en se ramenant au cas o\`u  $E$ est un module libre de rang deux, et en choisissant une base $\{e_{1}, e_{2}\}$ de $E$ : relativement \`a la base  $\{e_{1}^2, e_{1}e_{2}, e_{2}^2\}$  de ${\rm Sym}^2(E)$ et \`a la base $e_{1}\wedge e_{2} \otimes e_{1}\wedge e_{2}$ de 
$(\W^2E)^{\otimes 2}$, la matrice de cette application est 
$$
\begin{pmatrix}
0&0&2\\
0&- 1&0\\
2&0&0
\end{pmatrix}\eqno{\Box}
$$
\bigskip

\n {\bf 6.2.}\; Ainsi, \`a une section $\gamma$ de ${\mathcal Sym}^2(E) \otimes_{\OO}M$, cet isomorphisme associe une  forme $\varphi : {\mathcal Sym}^2(E)  \longrightarrow  L$  \`a valeurs dans le module inversible $L = (\W^2E)^{\otimes 2}\otimes M.$
\medskip

Si $\gamma$ est partout non nulle alors $\varphi$ est surjective, et r\'eciproquement. L'hypoth\`ese signifie, en effet,  que pour tout point $x \in X$, de corps r\'esiduel $\kappa(x)$, l'application 
$$
\gamma \otimes \kappa(x) : \kappa(x) \rr \sy^2(E\otimes \kappa(x)) \otimes_{\kappa(x)}(M\otimes \kappa(x))
$$
 est injective ; il revient au m\^eme de supposer que $\gamma$ est injective et que ${\rm Coker}(\gamma)$ est un $\OO$-module localement libre ([Bour] AC II, \S3.2, Prop.6, ou [EGA I] $0_{I} 6.7.4$). Mais c'est aussi \'evident si on utilise le lemme qui pr\'ec\`ede : cette dualit\'e montre que la non nullit\'e en chaque point de $\gamma$ et de $\varphi$ sont \'equivalentes, et, pour $\varphi$ elle est clairement \'equivalente \`a sa surjectivit\'e.
 \bigskip

\noindent {\bf  6.3. Proposition }\; {\it Soit } $
\gamma : \OO \, \longrightarrow \, {\mathcal Sym}^2(E) \otimes_{\OO}M
$ \; {\it une section partout non nulle, soit} $\varphi : {\mathcal Sym}^2(E)  \longrightarrow  L$ {\it la forme sym\'etrique ``duale'' de $\gamma$, et } $\A$ {\it  la $\OO$-alg\`ebre associ\'ee \`a $\varphi$, de sorte que $E$ est muni d'une structure de $\A$-module inversible. Alors, la suite}
$$
0 \longrightarrow \OO \stackrel{\gamma}{\longrightarrow} {\mathcal Sym}^2(E) \otimes_{\OO}M \stackrel{{\rm can.}}{\longrightarrow} {\mathcal Sym}_{\A}^2(E) \otimes_{\OO}M \longrightarrow 0
$$
{\it est exacte.}
\medskip

Ici encore, le plus simple est une v\'erification locale. On se place donc sur un ouvert affine o\`u $E$, $M$ et $L$ sont libres. On choisit une base $\{e_{1}, e_{2}\}$ de $E$, et une base $\beta$ de $M$ ; l'\'el\'ement $\varepsilon = (e_{1}\wedge e_{2})^{\otimes 2}\otimes \beta \, \in \, (\W^2E)^{\otimes 2}\otimes M = L$ est alors une base de ce module. La section $\gamma \in \sy^2(E)\otimes M$ s'\'ecrit 
$$
\gamma = (ae_{1}^2 + be_{1}e_{2}+ ce_{2}^2)\otimes \beta.
$$
On v\'erifie imm\'ediatement que la forme sym\'etrique associ\'ee \`a $\gamma$ est
$$
\varphi = 2a \,\varphi_{22} - b\, \varphi_{12} + 2c\, \varphi_{11}\leqno{(6.3.1)}
$$
o\`u $\varphi_{ij} : \sy^2(E) \longrightarrow L$ d\'esigne l'application donn\'ee par $\varphi_{ij}(e_{k}e_{l}) = \varepsilon$ si $\{i, j\} = \{k, l\}$, et $= 0$ sinon.
\medskip

\noindent  Reprenons  les notations de {\bf 4.2.}. L'alg\`ebre $\A$ associ\'ee \`a $\varphi$ est de la forme $\OO \oplus N$, o\`u le module $N$ est isomorphe \`a ${\mathcal Hom}(L, \W^2E)$, et  admet donc pour base l'\'el\'ement $\alpha$ d\'efini par $\alpha(\varepsilon) = e_{1}\wedge e_{2}$. La structure de $\A$-\`module sur $E$ est d\'etermin\'ee par l'action de $\alpha$, c'est-\`a-dire par un endomorphisme $\alpha_{E} : E \rr E$, qui est explicit\'e en {\bf 4.2.} ; on trouve
$$
\alpha_{E}(e_{1}) \; = \;  \frac{1}{2}b\, e_{1}\, + \, c\,e_{2}\, , \hspace{1,5cm}
\alpha_{E}(e_{2}) \; =\; - a\,e_{1}\, - \, \frac{1}{2}b\, e_{2}.\leqno{(6.3.2)}
$$

Par ailleurs, le noyau de l'application canonique $ {\mathcal Sym}^2(E) \otimes_{\OO}M \longrightarrow {\mathcal Sym}_{\A}^2(E) \otimes_{\OO}M$  est le $\OO$-module engendr\'e par les \'el\'ements de la forme $[ x \alpha_{E}( y) - \alpha_{E}(x)y]\otimes \beta  \in Sym^2(E)\otimes M$ ; on peut se limiter aux \'el\'ements $x, y$ faisant partie d'une base de $E$ ; ce noyau est donc engendr\'e par l'unique \'el\'ement
$$
[ e_{1}\alpha_{E}(e_{2}) - \alpha_{E}(e_{1})e_{2}]\otimes \beta =  \big[ e_{1}\,(- a\,e_{1}  - \frac{1}{2}b\, e_{2}) - ( \frac{1}{2}b\, e_{1}  +c\,e_{2})\,e_{2}\big ] \otimes \beta = - \gamma \hspace{2cm} \Box
$$ 
\vspace{7mm}

%%%%%%%%%%%%%%%%%%%%%
\noindent {\bf 7.\; Diviseurs de degr\'e deux sur les  fibr\'es en droites}
\bigskip

Soit $X$ un sch\'ema, $E$ un $\OO$-module localement libre de rang deux, et $M$ un $\OO$-module inversible. On consid\`ere une section partout non nulle
$$
\gamma : \OO \quad  \longrightarrow \quad {\mathcal Sym}^2(E) \otimes_{\OO}M.
$$

On peut voir $\gamma$ comme une famille, ind\'ex\'ee par $X$,  de polyn\^omes homog\`enes  de degr\'e 2 en deux ind\'etermin\'ees. Soit $P = {\bf P}(E)$ le fibr\'e projectif (en doites) associ\'e \`a $E$, et $p : P \rightarrow X$ le morphisme canonique. La section $\gamma$ d\'etermine  un diviseur effectif $D \subset P$, qui est fini, plat sur $X$, et localement de rang deux ; en particulier, le morphisme $D \rightarrow X$ est affine, et $D$ est donc le spectre de la $\OO$-alg\`ebre finie localement libre de rang deux $p_{\star}({\mathcal O}_{D})$. On se propose de d\'ecrire cette alg\`ebre.\\

\noindent {\bf 7.1.}\; La m\'ethode classique utilise les propri\'et\'es des images directes, rappel\'ees ci-dessous, o\`u  ${\mathcal O}_{P}(1)$ d\'esigne le quotient inversible ``fondamental" de $p^{\star}(E)$, et $F$ un $\OO$-module localement libre (cf. [EGA III] 2.1.16,\, ou [Hart] p.253, ex.8.3 et 8.4 ). 
$$
p_{\star}({\mathcal O}_{P}(m)\otimes p^{\star}(F)) \, \simeq \, \left \{ \begin{array}{c}
0, \, \hspace{2,2cm} {\rm si} \quad m < 0 \\
\sy^m(E)\otimes F, \quad {\rm si} \quad m \geq 0\\
\end{array}
\right.\leqno{(7.1.1)}
$$

$$
{\sf R}^1p_{\star}({\mathcal O}_{P}(m)\otimes p^{\star}(F)) \simeq \left \{ \begin{array}{c}
0,\hspace{2,4cm} {\rm si} \quad m \geq -1\\
(\W^2 E)^{ -1} \otimes_{\OO} F, \quad {\rm si}\quad  m = -2\\
\end{array}
\right.\leqno{(7.1.2)}
$$
Le diviseur $D$ associ\'e \`a $\gamma$ est d\'efini par $D = {\rm div}(s)$, o\`u $s$ est d\'eduite, par adjonction, de   
$$
\gamma :  \OO   \rightarrow  {\mathcal Sym}^2(E) \otimes_{\OO}M \simeq p_{\star}({\mathcal O}_{P}(2)\otimes_{{\mathcal O}_{P}}p^{\star}M) ; 
$$ 
autrement dit, $s$ est l'application compos\'ee
$$
{\mathcal O}_{P} \; \stackrel{p^{\star}(\gamma)}{\longrightarrow} \; p^{\star}p_{\star}({\mathcal O}_{P}(2)) \otimes_{{\mathcal O}_{P}}p^{\star}(M) \; \stackrel{{\rm can.}}{\longrightarrow}\; {\mathcal O}_{P}(2)\otimes_{{\mathcal O}_{P}}p^{\star}M .\leqno{(7.1.3)}
$$
On en tire la suite exacte
$$
0 \longrightarrow \; {\mathcal O}_{P}(-2)\otimes_{{\mathcal O}_{P}}p^{\star}M^{-1} \; \longrightarrow {\mathcal O}_{P} \; \longrightarrow \; {\mathcal O}_{D}\; \longrightarrow 0 . \leqno{(7.1.4)}
$$
Par image directe, on obtient la suite exacte
$$
0 \longrightarrow \; \OO \; \longrightarrow \; p_{\star}({\mathcal O}_{D})\; \longrightarrow \; {\mathsf R}^1p_{\star}({\mathcal O}_{P}(-2)\otimes_{{\mathcal O}_{P}}p^{\star}M^{-1}) \; \longrightarrow 0 \leqno{(7.1.5)}
$$
 Le faisceau conoyau ${\mathsf R}^1p_{\star}({\mathcal O}_{P}(-2)\otimes_{{\mathcal O}_{P}}p^{\star}M^{-1}) $ est, d'apr\`es (7.1.2),  isomorphe au module inversible 
 $$ N = (M\otimes \W^2 E)^{ -1} = {\mathcal Hom}(\W^2E \otimes M, \OO)
 $$
 Cela montre d\'ej\`a que $p_{\star}({\mathcal O}_{D})$ est une $\OO$-alg\`ebre finie localement libre de rang 2.
De plus, la suite (7.1.5) est scind\'ee (par ${1 \over 2}{\rm Tr}$), si bien que $p_{\star}({\mathcal O}_{D})$ est isomorphe, comme $\OO$-module,  \`a $\OO \oplus N$.  Par ailleurs, la suite (7.1.4), tensoris\'ee par ${\mathcal O}_{P}(1)$, donne, par image directe, un isomorphisme
$$
E\; = \; p_{\star}({\mathcal O}_{P}(1)) \quad \widetilde{\longrightarrow}\quad p_{\star}({\mathcal O}_{D}(1)).
$$
On en d\'eduit un morphisme de $\OO$-alg\`ebres, qui s'av\`ere \^etre injectif,
$$
p_{\star}({\mathcal O}_{D}) = p_{\star}{\mathcal End}_{{\mathcal O}_{D}}({\mathcal O}_{D}(1)) \; \longrightarrow \; {\mathcal End}_{\OO}(E).
$$
Autrement dit, $E$ est muni d'une structure de $p_{\star}({\mathcal O}_{D})$-module (inversible). Mais tout cela ne donne pas facilement la structure multiplicative sur $\OO \oplus N$, car les isomorphismes utilis\'es font intervenir la dualit\'e et des calculs \`a la \v Cech, qu'il serait malais\'e de suivre pour d\'egager le produit. On va montrer en ({\bf 7.3}) comment d\'ecrire cette multiplication dans le cadre propos\'e dans les paragraphes pr\'ec\'edents.
\vspace{0,5cm}

\n {\bf 7.2.}\; C'est le lieu de citer un r\'esultat, devenu classique, d\^u \`a Schwarzenberger ([Schwar], Thm 3, p.629), et qui se d\'eduit tr\`es simplement de ce qui pr\'ec\`ede.
\medskip

\n {\it Soit $X$un sch\'ema projectif sur un corps alg\'ebriquement clos de caract\'eristique $\neq 2$. On suppose que $\dim(X) \leq 2$. Alors,  tout $\OO$-module localement libre $E$ de rang $2$ est l'image directe $f_{\star}(\E)$ d'un inversible $\E$ sur un rev\^etement double $f : Y \rr X$.}
\medskip

 En effet, $X$ \'etant projectif, il existe un faisceau inversible tr\`es ample $M$ sur $X$ tel que $\sy^2(E)\otimes M$ soit engendr\'e par ses sections ; comme c'est un  $\OO$-module localement libre de rang $3 > \dim(X)$, ce module poss\`ede une section $\gamma : \OO \rr \sy^2(E)\otimes M$ partout non nulle, comme il d\'ecoule du \emph{easy lemma of Serre} ([Mumf], p.148). L'application $\varphi : \sy^2(E) \rr L$, duale de $\gamma$ au sens de  {\bf 6.1}, permet alors, en suivant {\bf 4.3} et {\bf 4.4}, de d\'efinir un rev\^etement double $f : Y \rr X$, et un  $\oo_{Y}$-module inversible $\E$ tels que $f_{\star}(\E) = E$.
\medskip

On verra en {\bf 9.1}\, que la conclusion est fausse pour un module ind\'ecomposable sur $\pp_{n}$ d\`es que $ n\geq 3$.
\bigskip

\noindent {\bf 7.3.}\; Revenons  aux notations introduites en {\bf 7.1.}, et \`a la section, suppos\'ee partout non nulle, 
$$
\gamma : \OO \quad  \longrightarrow \quad {\mathcal Sym}^2(E) \otimes_{\OO}M .
$$
Montrons comment le \S 4 permet de d\'ecrire l'alg\`ebre $p_{\star}(\oo_{D})$.

\noindent Le lemme {\bf 6.1}  permet d'associer \`a $\gamma$, par dualit\'e,  une application lin\'eaire 
$$
\varphi : {\mathcal Sym}^2(E)  \longrightarrow  L
$$
  \`a valeurs dans le module inversible $L = (\W^2E)^{\otimes 2}\otimes M.$ Cette application est surjective puisque $\gamma$ est partout non nulle. D'apr\`es le  \S 4, on associe \`a $\varphi$ un rev\^etement de rang deux $f : Y \rightarrow X$.
  \medskip
  
\noindent {\it  On va montrer l'existence d'un morphisme canonique de sch\'emas sur $X$, \, 
  $
  j : Y \longrightarrow P = {\bf P}(E)
  $, 
qui induit un isomorphisme de $Y$ sur le diviseur $D$ d\'efini par $\gamma$}.
  \medskip
  
  $$
\xymatrix{
Y \ar[dr]_{f }\ar[rr]^j && \pp(E) \ar[dl]^p \\
&X&
 }
$$

\noindent D\'esignons par $\E$ le module inversible sur $Y$ associ\'e \`a $E$, c.f. {\bf 4.3.} et {\bf 4.4.}, de sorte qu'on a
$$
f_{\star}(\E) \; = \; E.
$$
Sur $Y$, l'application surjective   $f^{\star}(E) = f^{\star}f_{\star}(\E) \rightarrow \E$ d\'etermine un morphisme
$j : Y \rightarrow P$ de sch\'emas sur $X$, caract\'eris\'e par l'isomorphisme 
$$
j^{\star}({\mathcal O}_{P}(1)) \simeq \E
$$
entre ${\mathcal O}_{Y}$-modules inversibles quotients de $f^{\star}(E)$. 
D'apr\`es [EGA II] 4.4.4 et 5.1.6, le morphisme $j$ est une immersion ferm\'ee puisque le morphisme $f$ est affine.
\medskip

Montrons  que  $j$ se factorise par le diviseur effectif  $D \subset P$ d\'efini par $\gamma$, c'est-\`a-dire, en utilisant (7.1.3),  que l'application compos\'ee suivante est nulle :

$$
{\mathcal O}_{Y} \; \stackrel{f^{\star}(\gamma)}{\longrightarrow} \; f^{\star}\sy^2_{{\mathcal O}_{X}}(E) \otimes_{{\mathcal O}_{Y}}f^{\star}(M) \; \stackrel{\psi \otimes 1}{\longrightarrow}\; \E^{\otimes 2}\otimes_{{\mathcal O}_{Y}}f^{\star}M, \leqno{(7.3.1)}
$$
o\`u 
$$
\psi :  f^{\star}(\sy^2(E)) \; {\longrightarrow}\; \E^{\otimes 2}
$$
est l'application canonique. Notons, comme  au \S1,  $
f_{\star}(\mathcal{O}_{Y}) \; =\; \A. 
$

\noindent En termes de $\OO$-modules, cette application $\psi$ s'\'ecrit
$$
\A \otimes_{\OO} \sy^2_{\OO}(E) \longrightarrow \sy^2_{\A}(E).
$$
Or, la proposition {\bf 6.3.} \'etablit l'exactitude de la suite
$$
0 \longrightarrow \OO \stackrel{\gamma}{\longrightarrow} {\mathcal Sym}^2(E) \otimes_{\OO}M\; \stackrel{{\rm can.}}{\longrightarrow} \;{\mathcal Sym}_{\A}^2(E) \otimes_{\OO}M \longrightarrow 0
$$
\noindent Il est alors clair que l'application (7.3.1) est nulle.

Finalement, comme les sch\'emas $Y$ et $D$ sont localement libres sur $X$, et de m\^eme rang, et que $j : Y \rightarrow D$ est une immersion ferm\'ee, on voit que $j$ est un isomorphisme.
\vspace{7mm}

\n {\bf 7.4.}\; A titre d'exemple, on va d\'eterminer l'alg\`ebre $p_{\star}(\oo_{D})$ dans le cas o\`u $X ={\rm Spec}(R)$ est affine, et o\`u $E = R^2$ et $M = R$. Le polyn\^ome homog\`ene consid\'er\'e s'\'ecrit alors
$$
\gamma \; = \; aX^2 + 2bXY + cY^2.
$$
Supposer que $\gamma$ est partout (sur ${\rm Spec}(R)$) non nulle revient \`a supposer que l'id\'eal $aR + bR + cR$ est \'egal \`a $R$.
D'apr\`es {\bf 6.1.} et les calculs utilis\'es dans la d\'emonstrations de  {\bf 6.3.}, on voit que la forme lin\'eaire $\varphi : \sy^2(E) \rr R$ associ\'ee \`a $\gamma$  s'\'ecrit
$$
\varphi = 2(c\varphi_{11} - b \varphi_{12} + a \varphi_{22}).
$$
L'alg\`ebre  $p_{\star}(\oo_{D})$  s'identifie \`a la sous-alg\`ebre de ${\rm End}_{R}(E) = {\bf M}_{2}(R)$ engendr\'ee par l'endomorphisme $\alpha_{E}$ de {\bf 4.2.}, soit, ici,  par la matrice 
$$
\begin{pmatrix}
b & -a\\
c&-b
\end{pmatrix}
$$
On trouve finalement que $p_{\star}(\oo_{D})$ est isomorphe \`a 
$$
A = R[T]/(T^2 - (b^2-ac)).
$$
On peut v\'erifier directement qu'on a un isomorphisme de sch\'emas
$$
{\rm Spec}(R[T]/(T^2-(b^2-ac))) \quad \widetilde{\longrightarrow}\quad {\rm Proj}(R[X, Y]/(aX^2+2bXY+cY^2)) \leqno{(7.4.1)}
$$
Il faut d'abord s'assurer que le $A$-module $E = R^2$ est inversible, puis que les \'el\'ements $\{e_{1}, e_{2}\}$ de la base canonique v\'erifient dans $E\otimes_{A}E$ la relation 
$$
a\, e_{1}\otimes e_{1} + 2b\, e_{1}\otimes e_{2} + c \,e_{2}\otimes e_{2} = 0.\leqno{(7.4.2)}
$$ 
V\'erifions ces deux points par un calcul direct. Pour le premier, il s'agit de montrer qu'il existe localement un \'el\'ement $z = (x, y) \in E$ tel que l'ensemble $\{z, \alpha_{E}(z)\}$ soit une base de $E$ comme $R$-module ; cela se traduit par l'inversibilit\'e (dans $R$) de
$$
{\rm det}\begin{pmatrix}x & bx-ay\\ y& cx - by
\end{pmatrix} \; = \; cx^2 - 2bxy +   ay^2.
$$
\noindent Consid\'erons donc le polyn\^ome $G(X, Y)= cX^2 -2bXY +aY^2$. 
Si l'\'el\'ement $a = G(0, 1)$ (resp. $c = G(1, 0)$) est inversible (dans $R$) alors l'\'el\'ement $(0, 1)$ (resp. $(1, 0)$) est un g\'en\'erateur de $E$ ; si $a-2b+c$ est inversible, alors $(1, 1)$ est un g\'en\'erateur. Comme l'hypoth\`ese sur les coefficients implique que $aR + (a-2b+c)R + cR = R$, on a bien v\'erifi\'e que le $A$-module $E$ est localement isomorphe \`a $A$.

\n Pour v\'erifier (7.4.2), il suffit de remarquer que le membre de gauche s'\'ecrit aussi
$$
- e_{1}\otimes(-a e_{1} -be_{2}) + (b e_{1} + c e_{2})\otimes e_{2}  = - e_{1}\otimes \alpha_{E}(e_{2}) + \alpha_{E}(e_{1})\otimes e_{2} ;
$$
cet \'el\'ement est nul puisque le produit tensoriel est effectu\'e sur $A$.

\newpage

%%%%%%%%%%%%%%%%%%%%%%%%%%%%%%

\noindent {\bf 8. \; Groupes de Picard d'un rev\^etement de l'espace projectif}
\bigskip

\noindent {\bf 8.1. Proposition}\; {\it Soit $X = {\bf P}_{n}$ l'espace projectif de dimension $n$ sur un corps de caract\'eristique $\neq 2$, et $f : Y \rightarrow X$ un rev\^etement double. Si $n \geq 3$, alors l'homomorphisme}
$$
f^{\star} : {\rm Pic}(X) \longrightarrow {\rm Pic}(Y)
$$
{\it est un isomorphisme.}
\bigskip

\noindent {\bf 8.2. Remarques}\; 1) \, Lorsque $n \leq 2$, on n'a plus un tel isomorphisme comme le montre, par exemple, le rev\^etement ${\bf P}_{1} \times {\bf P}_{1}\; \longrightarrow \; {\bf P}_{2}$.

\noindent 2) Sur le corps ${\bf C}$, et lorsque $Y$ est suppos\'e lisse, ce r\'esultat a \'et\'e d\'emontr\'e par Robert Lazarsfeld [Laz].

\noindent 3) La d\'emonstration propos\'ee ci-dessous utilise la th\'eorie des sections hyperplanes de Lefschetz, \'elabor\'ee par Grothendieck dans   [SGA 2], XI et XII  ; ce texte contient l'\'enonc\'e suivant (XII, Cor. 3.7), qui est tr\`es proche de {\bf 8.1.} : {\it si $Y$ est d'intersection compl\`ete (globale) dans un ${\bf P}_{r}$, alors ${\rm Pic}(Y)$ est libre engendr\'e par la classe de $\oo_{Y}(1)$} (aucune hypoth\`ese de lissit\'e sur $Y$). Je ne vois pas comment passer directement de cet \'enonc\'e \`a celui de {\bf 8.1}.
\bigskip

La d\'emonstration de {\bf 8.1} utilise plusieurs carr\'es cart\'esiens de la forme
$$
\begin{CD}
Y' @>f'>> X'\\
@VVV  @VVV\\
Y @>>f> X ,
\end{CD}\leqno{(\star)}
$$
o\`u $f$ (et donc $f'$) sont des rev\^etements doubles, et o\`u $X'$ est un diviseur effectif sur $X$. Cela conduit \`a un carr\'e commutatif
$$
\begin{CD}
{\rm Pic}(Y') @<f'^{\star}<< {\rm Pic}(X')\\
@AAA  @AAA\\
{\rm Pic}(Y) @<<f^{\star}< {\rm Pic}(X) ,
\end{CD}\leqno{(\star \star)}
$$
On va pouvoir comparer $f^{\star}$ et $f'^{\star}$ gr\^ace \`a des informations sur les fl\`eches verticales, fournies par le \og th\'eor\`eme de Lefschetz\fg\;  sous la forme donn\'ee par Grothendieck ([SGA 2] p.121), et dont voici l'\'enonc\'e.
\medskip

\noindent {\bf 8.3 \; Th\'eor\`eme}\; {\it Soit $X$ un sch\'ema alg\'ebrique projectif muni d'un module inversible ample $\LL$. On suppose que $X$ est de profondeur $\geq 2$ en ses points ferm\'es. Soit $X'$ le support d'un diviseur effectif d\'efini par une section de $\LL$.}
\medskip

\noindent {\it  \; Si ${\rm H}^1(X', \LL_{X'}^{\otimes - m}) = {\rm H}^2(X', \LL_{X'}^{\otimes - m}) = 0$ pour $m > 0$, et si les anneaux locaux des points ferm\'es de  $X$ sont parafactoriels, alors ${\rm Pic}(X) \longrightarrow {\rm Pic}(X')$ est bijectif.}
\bigskip

\n {\bf 8.4}\; En fait, dans la d\'emonstration de {\bf 8.1}, on disposera de conditions un peu plus fortes, et qui seront relatives \`a $X$, et non \`a son sous-sch\'ema $X'$. Gardons les hypoth\`eses g\'en\'erales de {\bf 8.3}.
\medskip

\noindent {\it  \;  Si ${\rm H}^j(X, \LL^{\otimes - m}) = 0$ pour $m > 0$, pour $j = 1, 2$ ou $3$, et si les anneaux locaux des points ferm\'es de  $X$ sont parafactoriels, alors ${\rm Pic}(X) \longrightarrow {\rm Pic}(X')$ est bijectif.}
\medskip

La suite exacte de faisceaux sur $X$
$$
0 \rr \LL^{-1} \rr \OO \rr \oo_{X'} \rr 0
$$
montre imm\'ediatement que les hypoth\`eses faites en {\bf 8.4} entra\^inent celles en {\bf 8.3}.
\bigskip

\noindent {\bf 8.5.}\; Pour m\'emoire :
\medskip

${\rm H}^j({\bf P}_{n}, \oo(r)) = 0$, \; pour $r \in {\bf Z}$ et $0 < j < n$.

\vspace{1cm}

\noindent {\bf 8.6.\; Lemme }\; {\it Soit $f : Y \rr X$ un rev\^etement double dont la multiplication $\mu : N^{\otimes 2} \rightarrow \OO$  est nulle, \emph{i.e.} dont le module inversible associ\'e $N$ est id\'eal de carr\'e nul dans $\oo_{Y}$. On suppose que}$$
{\rm H}^1(X, N) = {\rm H}^2(X, N) = 0.
$$
{\it Alors l'homomorphisme}
$
f^{\star} : {\rm Pic}(X) \longrightarrow {\rm Pic}(Y)
$
{\it est un isomorphisme.}
\medskip

\noindent Dans cet \'enonc\'e (classique), $X$ est un sch\'ema quelconque. La conclusion provient de la suite exacte
$$
0 \rr N \xrightarrow{x\mapsto 1+x} \oo_{X}^{\times} \rr \oo_{Y}^{\times} \rr 1
$$
\vspace{1cm}

\noindent {\bf 8.7} \; D\'emonstration de la proposition lorsque $n \geq 4$.

\noindent Soit $ f : Y \rightarrow X={\bf P}_{n}$ un rev\^etement double. Reprenons  les notations du \S  1 : $f_{\star}(\mathcal{O}_{Y}) = \oo_{X} \oplus N$. Le sch\'ema $X$ \'etant int\`egre, la multiplication $\mu : N^{\otimes 2} \rightarrow \OO$ est, soit nulle, soit injective. Si $\mu = 0$, $N$ est un id\'eal de carr\'e nul dans $\oo_{Y}$, et  le lemme pr\'ec\'edent  permet de conclure. Sinon  le lieu de diramation $X'$ est le diviseur  d\'efini  par l'id\'eal $\I = {\rm Im}(\mu) \subset \oo_{X}$.

Soit $Y' = f^{-1}(X') \stackrel{f'}{\rr} X'$ le rev\^etement induit au-dessus du diviseur $X'$, de sorte qu'on a le carr\'e cart\'esien suivant :

$$
\begin{CD}
Y' @>f'>> X'\\
@VVV  @VVV\\
Y @>>f> X ,
\end{CD}
$$

\noindent Montrons que ${\rm Pic}(X') \longrightarrow {\rm Pic}(Y')$ est un isomorphisme. 

\noindent Dans la d\'ecomposition  $f_{\star}(\mathcal{O}_{Y'}) = \oo_{X'} \oplus N'$, on a $N' = N/\I N$, et $N'$ est un id\'eal de carr\'e nul, si bien qu'on peut utiliser {\bf 8.6}, \`a condition d'avoir v\'erifi\'e que ${\rm H}^1(X', N') = {\rm H}^2(X', N') = 0$. Par d\'efinition, $N'$ s'ins\`ere dans la suite exacte
$$
0 \longrightarrow N\otimes N^{\otimes 2} \stackrel{1 \otimes \mu}{\longrightarrow} N \longrightarrow N' \longrightarrow 0
$$
Comme $X={\bf P}_{n}$, le module inversible $N$ est de la forme $N = \OO(-s)$, pour un entier $s > 0$, d'o\`u un isomorphisme $N\otimes N^{\otimes 2} \simeq \OO(-3s)$ ; enfin, comme $n \geq 4$, on a pour tout entier positif $r$, ${\rm H}^j(X, \OO(-r)) = 0$ pour $j = 1, 2, 3$. D'o\`u la nullit\'e de ${\rm H}^1(X', N')$ et de  $ {\rm H}^2(X', N')$.
\medskip

Consid\'erons maintenant les fl\`eches verticales, et v\'erifions les conditions de {\bf 8.4}. Comme $X$ est r\'egulier les anneaux locaux de ses points ferm\'es sont parafactoriels ; il en est de m\^eme pour $Y$ puisqu'il est de dimension $\geq 4$ et localement intersection compl\`ete ([SGA] 2, XI, 3.13, p.105).

Le faisceau ample $\LL$ de l'\'enonc\'e {\bf 8.4} est ici $\I^{-1} \simeq \oo_{X}(2s)$, donc les conditions d'annulation sont cons\'equences de {\bf 8.5}.

Ainsi, l'homomorphisme de restriction ${\rm Pic}(X) \rightarrow {\rm Pic}(X')$ est un isomorphisme.

Passons \`a $Y' \rightarrow Y$. Le faiceau ample $\LL$  \`a consid\'erer ici est  $\I \oplus \I N \simeq  \oo_{X}(-2s) \oplus \oo_{X}(-3s)$, et il est clair que l'on a les m\^emes propri\'et\'es d'annulation ; elles entra\^inent que ${\rm Pic}(Y) \rightarrow {\rm Pic}(Y')$ est un isomorphisme.

Finalement, le carr\'e ($\star \star$) montre que ${\rm Pic}(X) \rightarrow {\rm Pic}(Y)$ est un isomorphisme.
\vspace{1cm}

\noindent {\bf 8.8}\; On suppose maintenant que $X = {\bf P}_{3}$.

\noindent On r\'ealise $X$ comme un hyperplan dans $Q = {\bf P}_{4}$. Le rev\^etement donn\'e $f : Y \rightarrow X$ est d\'etermin\'e par l'application injective $\mu : \oo_{X}(-2s) \rightarrow \oo_{X}$, c'est-\`a-dire par une section de $\oo_{X}(2s)$. La suite exacte
$$
0 \longrightarrow \oo_{Q}(2s-1) \rr \oo_{Q}(2s) \rr \oo_{X}(2s) \rr 0
$$
et les relations {\bf 8.5}, montrent  que $\mu$ se rel\`eve (sans unicit\'e) en une application
$$
\mu_{Q} : \oo_{Q}(-2s) \rr \oo_{Q}.
$$
Cette application permet de d\'efinir un rev\^etement double $g : Z \rr Q$ dont la restriction \`a l'hyperplan $X \subset Q$ est $f$. On a donc le carr\'e cart\'esien suivant :
$$
\begin{CD}
Y@>f>> X\\
@VVV  @VVV\\
Z @>>g> Q 
\end{CD}
$$
On en d\'eduit le carr\'e commutatif
$$
\begin{CD}
{\rm Pic}(Y) @<f^{\star}<< {\rm Pic}(X)\\
@AAA  @AAA\\
{\rm Pic}(Z) @<<g^{\star}< {\rm Pic}(Q) ,
\end{CD}
$$
Comme $Q = {\bf P}_{4}$, la premi\`ere partie ({\bf 8.7}) montre que $g^{\star}$ est un isomorphisme. Comme $X$ est un hyperplan dans $Q$, on dispose de la suite exacte
$$
0\rr \oo_{Q}(-1) \rr \oo_{Q} \rr \OO \rr 0
$$
Par suite les conditions de {\bf 8.4} sont satisfaites, et  ${\rm Pic}(Q) \rr {\rm Pic}(X)$ est donc un isomorphisme. 

\noindent Comme l'id\'eal de $Y$ dans $Z$ est $\oo_{Q}(-1) \oplus \oo_{Q}(-1)\otimes N = \oo_{Q}(-1) \oplus \oo_{Q}(-1-s)$, les conditions d'annulation de {\bf 8.4} sont aussi satisfaites sur $Z$, donc
${\rm Pic}(Z) \rr {\rm Pic}(Y)$ est un isomorphisme. Cela permet de conclure.
\vspace{1,5cm}

%%%%%%%%%%%%%%%%%%%%%%%%%%%

\noindent {\bf 9. \; Application aux fibr\'es de rang deux sur les espaces projectifs.}
\bigskip

\noindent {\bf 9.1. \; Proposition}\; {\it Soit $X = {\bf P}_{n}$ l'espace projectif de dimension $n \geq 2$ sur un corps alg\'ebriquement clos de caract\'eristique $\neq 2$. Soit $E$ un fibr\'e de rang 2 sur $X$, et $\varphi : \sy^2(E) \rightarrow \OO(r)$ une application r\'eguli\`ere.

\noindent {\it i})\; Si $E$ est ind\'ecomposable, alors $ r > c_{1}$.

\noindent {\it ii})\; Si $\varphi$ est surjective et si $n \geq 3$, alors $E$ est d\'ecomposable.}
\bigskip

\noindent  L'hypoth\`ese sur $\varphi$ signifie ici qu'en tout point $x \in X$ tel que $\dim(\oo_{X, x}) < {\rm rang}(\sy^2(E)) = 3$, l'application $\oo_{X, x}$-lin\'eaire  $\varphi_{x}$ est surjective.

\n Comme d'habitude, l'entier $c_{1}$ est d\'efini par $\W^2 E = \OO (c_{1})$.
\medskip

\noindent Suivant le \S 4.1, on introduit l'application associ\'ee \`a $\varphi$
$$
u : N\otimes E \, \longrightarrow \, E
$$
o\`u $N = \OO (c_{1}-r)$ ; rappelons que $1_{\kappa(x)}\otimes u$ est non nulle en tout point $x \in X$ o\`u $\varphi$ est non nulle, {\it i.e} o\`u $\varphi_{x}$ est surjective.
\bigskip

\noindent Montrons que si $u$ n'est pas injective, alors $E$ est d\'ecomposable.
\medskip

\noindent Supposons donc que le module  $M = {\rm Ker}(u)$ ne soit pas nul, et montrons qu'alors il est inversible. L'exactitude de la suite $0\rightarrow M \rightarrow N\otimes E \stackrel{u}{\rightarrow} E$ montre d'abord que $M$ est est g\'en\'eriquement de rang 1 puisque $u$, comme $\varphi$, est g\'en\'eriquement non nulle ; comme $X$ est lisse, le dual $M^{\vee}  = {\mathcal Hom}(M, \OO)$ est un $\OO$-module inversible.

\begin{quotation}
Plut\^ot qu'un liste de r\'ef\'erences \'eparpill\'ees qui conduiraient \`a ce r\'esultat bien connu, voici un d\'emonstration directe : soit $R$ un anneau local r\'egulier, de corps des fractions $K$ et $M$ un $R$ module de type fini tel que $K\otimes_{R} M$ soit de rang 1. Alors $M^{\vee} = {\rm Hom}(M, R)$ est libre. Soit, en effet, $v : M \rightarrow R$ une forme lin\'eaire telle que l'id\'eal $v(M)$ soit maximal parmi les id\'eaux de ce type ; montrons que $v$ est une base du dual $M^{\vee}$. Notons d'abord que, pour un \'el\'ement non nul $t \in R$, la relation $v(M) \subset tR$ implique que $t$ est inversible ; en effet, elle entra\^ine l'existence de $v' \in M^{\vee}$ tel que $v = tv'$ ; comme $v(M) = tv'(M) \subset v'(M)$, la maximalit\'e de $v(M)$ implique que $tv'(M) = v'(M)$, donc que $t$ est inversible. Soit $w \in M^{\vee}$, et $\xi \in K$ l'\'el\'ement tel que $w = \xi\, v$ ; il faut montrer que $\xi$ est dans $R$, c'est-\`a-dire, $R$ \'etant normal, que $\xi \in R_{\mathfrak{p}}$ pour tout id\'eal premier $\mathfrak{p}$ de hauteur 1 ; mais, $R$ \'etant factoriel, un tel id\'eal est principal  : $\mathfrak{p} = tR$ ; ce qui pr\'ec\`ede entra\^ine que  $v(M) \nsubseteq \mathfrak{p}$, donc que $R_{\mathfrak{p}} =  v(M)_{\mathfrak{p}}$ ; par suite, on a $\xi R_{\mathfrak{p}} = \xi v(M)_{\mathfrak{p}} = w(M)_{\mathfrak{p}} \subset R_{\mathfrak{p}}$ ; d'o\`u le r\'esultat.
\end{quotation}

La suite exacte qui d\'efinit $M$ montre aussi que ce module est reflexif ; c'est donc un $\OO$-module inversible.

Montrons que $ {\rm Im}(u)$ est un module inversible. Comme $X$ est int\`egre, et que $u$ n'est pas injectif, on a $\W^2u = 0$ ([Bour], A III 8.2, Prop. 3) ; par suite, la multiplication $\mu : N^{\otimes 2} \rightarrow \OO$ est nulle (4.1.4), et on a $u\circ (1_{N}\otimes u) = 0$ ; on en tire que l'application $1_{N}\otimes u : N^{\otimes 2} \otimes E \rr N\otimes E$ se factorise en 
$$
N^{\otimes 2} \otimes E \longrightarrow N\otimes  {\rm Im}(u)\, \subset \, M \subset N\otimes E.
$$
Posons $M' =  {\rm Im}(u)$. Comme l'application  $1_{N}\otimes u$ est non nulle en les points de codimension $\leq 2$ (tout comme $\varphi$), il en est de m\^eme de l'inclusion $N\otimes M' \subset M$ ; mais $M$ \'etant localement isomorphe \`a $\OO$, Nakayama nous enseigne que la non nullit\'e de $N\otimes M'\otimes \kappa(x) \rightarrow M\otimes \kappa(x)$ \'equivaut \`a la surjectivit\'e de $N\otimes M'_{x} \rightarrow M_{x}$ ; par suite le support du quotient $Q = M / (N\otimes M')$ est de codimension $\geq 3$ ; mais, par ailleurs, si le module $Q$ n'est pas nul, sa dimension projective est $\leq 2$ puisqu'on a la suite exacte
$$
0\rightarrow N\otimes M \longrightarrow N^{\otimes 2}\otimes E\,  \stackrel{1\otimes u}{\longrightarrow}\, M \rightarrow Q \rightarrow 0 .
$$
Ceci est impossible sur un sch\'ema r\'egulier (Th\'eor\`eme de Auslander-Buchsbaum, [Bour], AC, X p.45) ; donc $N\otimes M'\, = \, M$, et par suite $M'$ est un module inversible ; mais alors $N\otimes E$ appara\^it comme une extension de $M'$ par $M$, qui sont des faisceaux inversibles :
$$
0\rightarrow M \rightarrow N\otimes E \stackrel{u}{\rightarrow} M' \rightarrow 0
$$
Une telle extension correspond \`a un \'el\'ement de ${\rm Ext}^1(M', M) = {\rm H}^1({\bf P}_{n}, N)$, et ce groupe est nul d\`es que $n \geq 2$. Cette extension est donc scind\'ee, ce qui montre que si $u$ n'est pas injective, alors $E$ est d\'ecomposable.\medskip

 \n Les hypoth\`eses impliquent donc que l'application $
u : N\otimes E \, \longrightarrow \, E
$ est injective ; mais alors $\wedge^2 u$ est elle aussi injective, ainsi que $\mu : N^{\otimes 2}= \OO(2(c_{1}-r)) \, \rightarrow \OO$ ; donc $r \geq c_{1}$.
\medskip

Il reste \`a \'ecarter le cas o\`u $r = c_{1}$. Soit $Z \subset X$ le sch\'ema des z\'eros de $\varphi$. Notons $f : Y \rr X$ le rev\^etement double associ\'e \`a $N$ et $\mu$. D'apr\`es {\bf 4.3}, il existe un faisceau coh\'erent $\E$ sur $Y$ tel que $f_{\star}(\E) = E$, et d'apr\`es {\bf 4.4}, ce faisceau est inversible sur l'ouvert $f^{-1}(X - Z)$.

\n Si $r =c_{1}$,  la multiplication $\mu$ est un isomorphisme, donc le rev\^etement $f : Y \rightarrow X$ est \'etale fini ; comme $X = {\bf P}_{n}$ est simplement connexe ([SGA] 1, p.219), on a  $Y = X \sqcup X$, et le faisceau $\E$ sur $Y$ est d\'ecompos\'e en $\E' \oplus \E''$ ; comme $\E$ est inversible sur l'ouvert  $f^{-1}(X - Z)$, chacun des composants $\E'$ et $\E''$ est non nul ; mais cela implique que $E$ est d\'ecompos\'e.
\bigskip

\noindent Montrons {\it ii}).

\noindent Notons encore $f : Y \rr X$ le rev\^etement associ\'e \`a $\varphi$, et $\E$ le faisceau sur $Y$ tel que $f_{\star}(\E) = E$. Si $\varphi$ est surjective, alors le faisceau $\E$ sur $Y$ est inversible ({\bf 4.4}). Comme $n \geq 3$, la proposition {\bf 8.1}  montre que l'homomorphisme ${\rm Pic}(X) \rr {\rm Pic}(Y)$ est un isomorphisme ; on a donc $\E = f^{\star}(L)$ pour un inversible  $L$  sur $X$ ; mais alors le module $E = f_{\star}f^{\star}(L) = L \oplus \, (L\otimes N)$ est d\'ecompos\'e.
\vspace{1cm}

%%%%%%%%%%%%%%%%%%%%%%%%%%%
{\bf Bibliographie}
\medskip

\noindent [Anan], S. ANANTHARAMAN, {\it Sch\'emas en groupes $\ldots$} Bull. Soc. Math. France, M\'emoire 33, (1976)

 \noindent [Bour] A,   N. BOURBAKI, Alg\`ebre 
 
\noindent [Bour] AC,   N. BOURBAKI, Alg\`ebre commutative

 \noindent [Fer 1], D. FERRAND, {\it Un foncteur norme} Bull. Soc.math. France, {\bf 126}, (1998) p.1-49
 
 \noindent [Fer 2],  D. FERRAND, {\it Conducteur, descente et pincement}, Bull.Soc.math.France, {\bf 131}, (2003), p.553-585

 \noindent [EGA],  A. GROTHENDIECK, r\'edig\'es  avec la collaboration de J. DIEUDONN\'E , \'El\'ements de G\'eom\'etrie alg\'ebrique, Springer-Verlag (1971), et Publ.Math.IHES, Paris (1960-1967)
 
 \noindent [SGA] 1,  A. GROTHENDIECK, Rev\^etements \'etales et groupe fondamental, {\it S\'eminaire du Bois-Marie 1961}, Documents Math\'ematiques 3, Soc.math. de France (2003)

 \noindent [SGA] 2,  A. GROTHENDIECK, Cohomologie locale des faisceaux coh\'erents \ldots , {\it S\'eminaire du Bois-Marie 1962}, Documents Math\'ematiques 4, Soc.Math. de France (2005)

 \noindent [Hart],   R. HARTSHORNE, Algebraic Geometry, Springer-Verlag (1977)
 
 \n [Knes], M. KNESER, {\it Composition of Binary Quadratic Forms}, J. of Number Thory, {\bf 15}, (1982), p.406-413
 
\noindent [Laz], R. LAZARSFELD, {\it A Barth-Type Theorem for Branched Coverings of Projective Space}, Math. Ann. {\bf 249}, (1980) p.153-162 

\n [Mumf], D. MUMFORD, {\it Lectures on Curves on an Algebraic Surface}, Ann. of Math. Studies, Number 59, (1966), Princeton Univ. Press

\noindent [Roby], N. ROBY, {\it Lois polyn\^omes et lois formelles en th\'eorie de modules}, Ann. scient. \' Ec. Norm. Sup., t.80, (1963), p.213-348.

\n [Schwar], R. L. E. SCHWARZENBERGER, {\it Vector Bundles on the Projective Plane}, Proc. London Math.Soc. {\bf 11}, (1961), p.623-40.

\end{document}